\documentclass[a4paper]{amsart} % [a4paper,10pt,openright]{article}
\usepackage{graphicx}
\usepackage[english]{babel}
\usepackage{multicol}
\usepackage{amsmath}
\usepackage{amssymb}
\usepackage{amsfonts}
\usepackage{mathrsfs}
\usepackage{graphicx}
\usepackage{fancyhdr}
\usepackage{newlfont}
\usepackage{indentfirst}
\usepackage{geometry}
\usepackage{enumerate}
\usepackage{color}
\usepackage[foot]{amsaddr}
\usepackage[font={small,it}]{caption}

\begin{document}

\title{From receptive profiles to a metric model of V1}
 \date{}
\author{Noemi Montobbio$^\ast$}
\address{$^\ast$Dipartimento di Matematica, Universit\`{a} di Bologna, Italy. noemi.montobbio2@unibo.it}
\author{Giovanna Citti$^\ddagger$}
\address{$^\ddagger$Dipartimento di Matematica, Universit\`{a} di Bologna, Italy. giovanna.citti@unibo.it}
\author{Alessandro Sarti$^\dagger$}
\address{$^\dagger$CAMS, CNRS - EHESS, Paris, France. alessandro.sarti@ehess.fr}

\begin{abstract}
In this work we show how to construct connectivity kernels induced by the receptive profiles of simple cells of the primary visual cortex (V1). These kernels are directly defined by the shape of such profiles: this provides a metric model for the functional architecture of V1, whose global geometry is determined by the reciprocal interactions between local elements. Our construction adapts to any bank of filters chosen to represent a set of receptive profiles, since it does not require any structure on the parameterization of the family.\\
The connectivity kernel that we define carries a geometrical structure consistent with the well-known properties of long-range horizontal connections in V1, and it is compatible with the perceptual rules synthesized by the concept of association field. These characteristics are still present when the kernel is constructed from a bank of filters arising from an unsupervised learning algorithm.
\end{abstract}

\maketitle

\section{Introduction}\label{intro}
 
 The primary visual cortex (V1) implements the first stage of cortical processing of the visual information, and it is the most studied and the best understood among the visual areas in the brain. An essential concept related to the structure of V1 is that of \emph{receptive profile} (RP), i.e. the function modeling the impulse response of a neuron. Indeed, the primary step of computation taking place in V1 involves a class of neurons, called \emph{simple cells}, which act approximately as linear filters on the retinal image. In other words, the action of such a neuron on a visual stimulus can be expressed, to a first approximation, as the convolution between the image and the RP of the cell.

 A classical model for the set of receptive profiles of V1 simple cells is provided by a bank of bidimensional Gabor filters \cite{jonpal,daug,lee}. Such a family can be indexed by a set of parameters showing a group structure: this fact is at the basis of a number of models, which identify the set of simple cells with such space of parameters and define on it a differential structure, designed to be compatible with the above-mentioned group law \cite{petitond,cs06,scp}. This idea can in principle be replicated as long as the bank of filters modeling the RPs can be parameterized by a group. Yet, this condition is not verified if one considers for example a bank of filters learned through automatic learning procedures (see e.g. \cite{sparse} and \cite{anselmipoggio}). It is nonetheless crucial for a cortex model to describe the \emph{functional architecture} ruling the geometry of intra-cortical connections: these are indeed thought to be responsible for the integration of local features into contours, shapes, objects. Specifically, the contextual influences modulating V1 neurons' activity have been described through the concept of \emph{association field} \cite{field}, which specifies the geometrical properties of the influences between neurons even with markedly separated RPs, based on their reciprocal position and orientation. The long-range, orientation specific horizontal connections which are known to take place in V1 are believed to be at the basis of such mechanisms of perceptual grouping \cite{bosking}. \\

 The aim of this paper is to show that a bank of receptive profiles can induce such a functional architecture, with or without the presence of a group structure.\\
 We propose to express the connectivity strength through a kernel defined onto the family of RPs and generated by the RPs themselves. Specifically, we consider a bank of filters $\{\psi_p\}_{p \in \mathcal{G}}$ $\subseteq L^2(\mathbb{R}^2)$, where $\mathcal{G}$ is any set of parameters indexing the family. As in \cite{cs06}, we identify each filter $\psi_p$ with the corresponding index $p \in \mathcal{G}$. We shall refer to $\mathcal{G}$ as the \emph{feature space}: intuitively, one may think that each element $p \in \mathcal{G}$ encodes the features extracted by the corresponding filter $\psi_p$ when it is applied to some image. We then define, for each couple of points $p,p_0 \in \mathcal{G}$, the \emph{generating kernel}
 $$K(p,p_0) := Re \left(\int_{\mathbb{R}^2} \psi_p(x,y) \: \overline{\psi_{p_0}(x,y)} \: dx \: dy \right),$$
 i.e. $K(p,p_0)$ is defined to be the real part of the $L^2$ scalar product between the filters $\psi_{p}$ and $\psi_{p_0}$. This kernel allows to define a natural metric structure on the feature space. We then construct, through an iterative mechanism of repeated integrations against $K$, a wider connectivity kernel from which we are also able to recover association-field-like patterns. On the other hand, in differential cortex models the spreading of neural activity across V1 is typically described by means of diffusion equations expressed through second order operators associated to the metric. In our work \cite{metric}, we consider a diffusion process defined in terms of the distance, and we study its relation with the iterative procedure described in this paper. The technique introduced in \cite{metric} allows to extend the diffusion-based approach used in differential models to the present metric setting, including possibly non-differential metrics, while still incorporating the sub-Riemannian model introduced in \cite{cs06} as a limit case.\\
 
 In the second part of the paper, we present some results and numerical simulations to demonstrate the applicability of our model to both the classical framework of a family of Gabor filters (including a spatiotemporal example), and the situation where the RPs of simple cells are represented by a bank of filters obtained through an unsupervised learning algorithm. The application of this construction to \emph{endstopped} RPs also allowed to recover the property of curvature selectivity observed \cite{dobbins87} in these neurons. Our results suggest that \emph{the geometry of the horizontal connections in V1 can indeed be obtained from the shape of the receptive profiles}.\\
 The promising outcomes obtained even in the case of a non-structured bank of learned filters led us to the idea of applying our construction in the context of convolutional neural networks (CNNs) \cite{lecun}. We thus propose, as a future development, the application of our model to insert lateral connections induced by the feedforward filters in the convolutional layers of a CNN.

\section{The functional architecture of V1}\label{functional}

In this section we first recall the basic notions of receptive field and receptive profile of a visual neuron, with a particular focus on V1. We then review the link between the mechanisms of perceptual grouping and the neural circuitry of horizontal connections in V1. Lastly, we briefly outline the main concepts underlying some existing cortex models, mainly based on differential structures.

\subsection{Receptive fields and receptive profiles}\label{rfrp}

 Every subcortical or cortical neuron is linked, through the connectivity of retino-geniculo-cortical pathways, to a local domain $D$ of the retina, called the \emph{receptive field} (RF) of the cell. The \emph{receptive profile} (RP) of a visual cell is a function $\psi(x,y)$ defined on $D$, measuring quantitatively the impulse response of the neuron.\\
Certain classes of cells act approximately as \emph{linear filters} on the optic signal, meaning that their receptive profile provides, up to a first approximation, the necessary information to measure their reaction to a general stimulus (expressed as a function $I(x,y)$ of the retinal coordinates): this is indeed given by a \emph{convolution}
$$ \int_D I(x,y)\psi(x-x_0,y-y_0)dxdy.$$
Here, the RP is expressed in local coordinates around $(x_0,y_0)$, the retinal location around which it is concentrated: that is, the response of the cell to a punctual stimulus at $(x_0,y_0)$ is given by $\psi(0,0)$ .

\subsection{The structure of V1}\label{v1}

 Let us now focus on the structure of the first of the visual cortical areas. V1 neurons can be divided into two main classes, as first discovered by Hubel and Wiesel in 1962 \cite{HW}: these are called \emph{simple} and \emph{complex} cells.\\
 Simple cells are the first neurons in the visual pathway which are sensitive to orientation.
 \begin{figure}[ht]
\centering
 \includegraphics[width=0.5\textwidth]{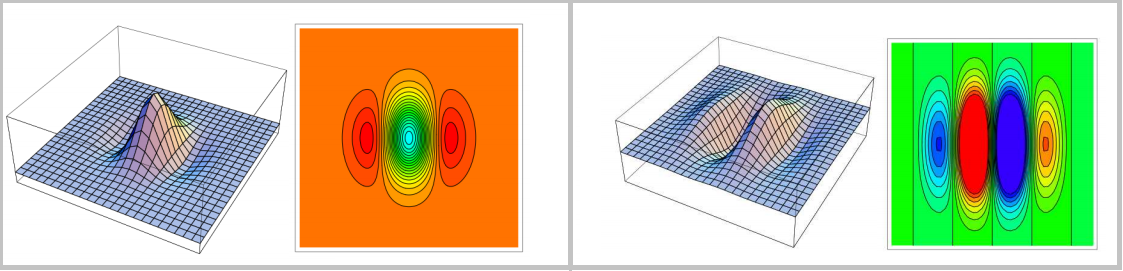}
 \caption{Two examples of RPs of simple cells of V1. Source: \cite{petitot}.}\label{RP_V1}
\end{figure}
The basic structure of their RPs consists of two to five tapering rows of alternating excitatory and inhibitory subregions \cite{daug}. Two examples of such profiles are shown in Figure \ref{RP_V1}. As already mentioned, the responses of simple cells display orientation and positional selectivity: each of them fires at an optimal orientation, giving progressively weaker responses as the orientation of the stimulus shifts sub-optimally. This behavior is classically described by modeling the RPs of this class of cells through a family of \emph{Gabor filters}, defined as complex exponential functions modulated by bidimensional Gaussian functions \cite{jonpal,daug,lee}. The filters vary in retinal position, orientation and scale in order to describe the whole set of profiles. Starting from the basic filter
\begin{equation}\label{mother} \psi_{0,0,0,1} (u,v) = \exp\left(\frac{2\pi i u}{\lambda}\right) \exp\bigg(-\frac{u^2+v^2}{2\sigma^2}\bigg),
 \end{equation}
 one can obtain the whole bank by translations, rotations and dilations. These are expressed by the following parameters: $(x,y) \in \mathbb{R}^2$ represents the position at which the filter is translated, $\theta \in S^1$ is the angle of the rotation and $\sigma>0$ is the factor of dilation. Each combination of these transformations determines a filter $\psi_{x,y,\theta,\sigma}$ indexed by such parameters.\\
 A different kind of processing is carried out by complex cells: unlike simple cells, their response cannot be explained through the linear action of a receptive profile on an image. These cells still display orientation sensitivity, but are \emph{phase invariant}, i.e. they respond equally to optimally oriented bars independent of where they are placed inside the RF.
 A number of findings (see \cite{complex} for a review) led many authors to model the response of a complex cell as a square sum of simple cells with similar orientation and scale but with phases that differed in 90 degrees. Such a couple of cells is commonly referred to as a \emph{quadrature pair}. This characterization of the development of phase invariance in complex cells is often called \emph{energy model}.
 
 However, the interactions between V1 neurons are not limited to the classical hierarchy of simple - to - complex cells.
 \begin{figure}[h]
\centering
 \includegraphics[width=0.3\textwidth]{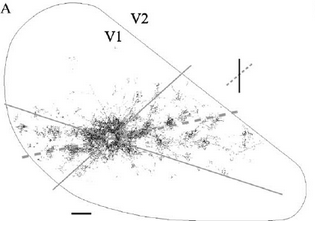}
 \caption{The horizontal connections departing from a V1 neuron extend along an axis according to the preferential orientation of the cell. Source: \cite{bosking}. }\label{bosking}
\end{figure}
 Indeed, there exist two main types of intracortical connectivity affecting the responses of V1 neurons. Short-range connections, linking neurons whose RPs have the same retinal position but different orientations, act by selecting the orientation causing the maximal response to a stimulus and suppressing the others. On the other hand, lateral (or ``horizontal'') connections associate neurons sharing the same preferred orientation, but linked to different positions on the retina. The qualitative and quantitative properties of the latter have been investigated in several neurophysiological experiments, the best known of which are the ones performed by W. Bosking et al. in 1997 \cite{bosking} (see Figure \ref{bosking}): these connections show a marked orientation specificity and spread for millimeters horizontally along the cortex, linking cells whose RF do not overlap. Such attributes make them a likely candidate to account for the processes of perceptual grouping carried out by our visual system.
 
 \subsection{Contextual influences and association fields}\label{context}

The human visual system is highly efficient in detecting and segregating objects. One of the main aspects in the perception of shapes is the identification of contours: this cannot be explained through the classical concept of RF \cite{gromin,neumin,angelucci03}, which refers to the extraction of \emph{local} image features from the visual scene. In other words, the activity of cells is not only influenced by inputs within their RF, but rather it strongly depends on the context. \\
\begin{figure}[h]
\centering
 \includegraphics[width=0.5\textwidth]{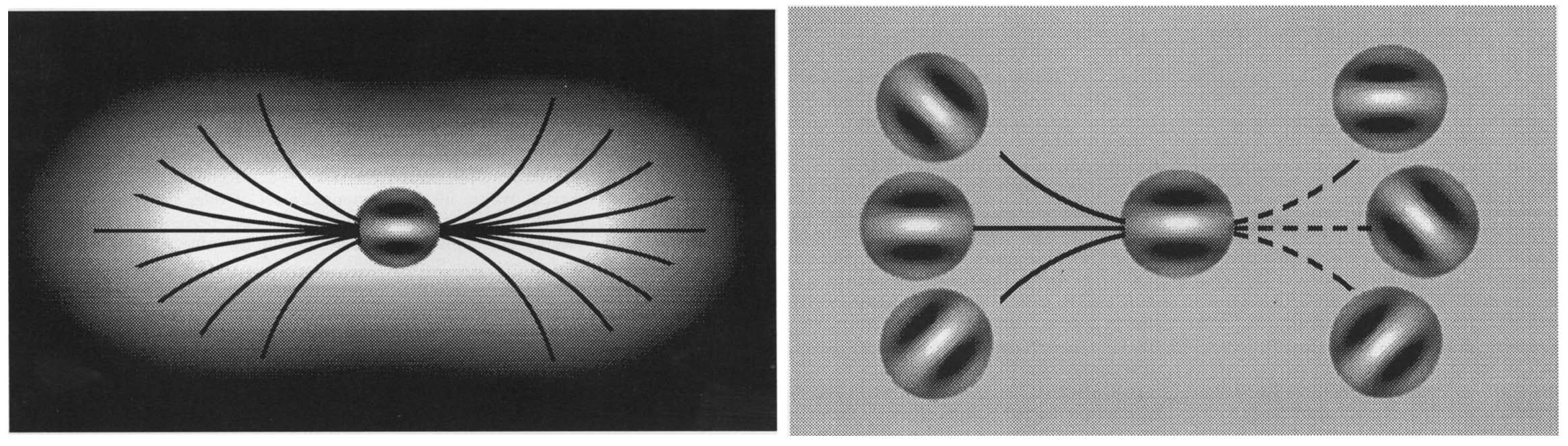}
 \caption{A schematic representation of the concept of association field introduced in \cite{field}. Left: the ``area of influence'' of a cell and some of the possible paths departing from it. Right: examples of neurons for which grouping with the central cell occurs (solid curves) or does not occur (dashed curves).}\label{ass_field}
\end{figure}
As far as V1 is concerned, it has been proposed (Field, Hayes and Hess, 1993 \cite{field}) that the interaction of local edge elements in this area can be described through an \emph{association field}, consisting of a region around each V1 cell where the activation of other neurons characterized by certain preferred orientations tends to be enhanced by the excitation of the central cell. See Figure \ref{ass_field}. A number of psychophysical experiments, described in \cite{field}, suggested that the area covered by this activation region around a neuron is considerably wider than the RF of the cell.\\
In V1, the interactions mediated by the association field show facilitatory influences for neurons whose RPs have a similar preferred orientation; however, there may exist an analogue entity in other cortical areas \cite{gilwu}, which links neurons sharing different and more abstract features. In other terms, each area could be endowed with an \emph{intrinsic} correlation kernel, closely linked to the type of analysis performed in it, yielding corresponding perceptual patterns.

 As just mentioned, the association field observed in V1 is believed to be intrinsic to this area: a possible anatomical framework is represented by the long-range horizontal connections \cite{petitond,gilbert}, as they fittingly link neurons with widely separated RPs whose reciprocal orientations are collinear or co-circular.

\subsection{Previous models}\label{previous}
 
 In this section, we give an overview on some previous mathematical models describing the connectivity of the visual cortex and the laws of perceptual organization.\\
 
 Neuromathematical models of V1 often adopt a differential geometry approach to modeling the functional architecture of this area. A breakthrough idea has been that of representing V1 as a fiber bundle with the space of retinal locations as a basis. This approach was originally proposed by Petitot and Tondut \cite{petitond}, with the definition of a left-invariant sub-Riemannian structure on the Heisenberg group. In 2006, Citti and Sarti wrote the model in the Lie group $SE(2) = \mathbb{R}^2 \times S^1$ by requiring the invariance under roto-translations \cite{cs06}. Such sub- Riemannian structure has been shown to be naturally induced by the action of a bank of Gabor filters on a visual stimulus. In particular, consider a Gabor function 
 \begin{equation*}
  \Psi(x,y,\theta)= \exp\left(-\frac{\xi^2+\nu^2}{\sigma^2}\right) \exp\left( i\frac{\nu}{\sigma^2} \right), 
  \quad \text{ with } \begin{cases}
        \xi & = x\cos\theta+y\sin\theta,\\
        \nu & = -x\sin\theta+y\cos\theta. \end{cases}
 \end{equation*}
 The outcome of filtering an image $I$ with $\Psi$ can be locally approximated as
 \begin{equation}
  O(x,y,\theta) = X_3 \exp\left(-\frac{\xi^2+\nu^2}{\sigma^2}\right) \ast I =: X_3 \: I_\sigma,
 \end{equation}
 where \begin{equation}
        X_3 = X_3(\theta) = -\sin\theta \partial_x + \cos\theta \partial_y.
       \end{equation}
 $I_\sigma$ is a smoothed version of $I$, obtained by convolving it with a Gaussian kernel.\\
 This procedure \emph{lifts} the 2D image domain to the 3D space $\mathbb{R}^2 \times S^1$: each point $(x,y) \in \mathbb{R}^2$ is sent to a point $(x,y,\overline{\theta}) \in \mathbb{R}^2 \times S^1$ such that $\overline{\theta}$ is a local maximum point of $\theta \mapsto O(x,y,\theta)$, so that the whole image domain is lifted to the set
 \[
  \{ (x,y,\overline{\theta}) \; : \; O(x,y,\overline{\theta}) = \max_\theta O(x,y,\theta) \} \; \subseteq \; \mathbb{R}^2 \times S^1.
 \]
 This ``non-maximal suppression'' principle is based on experimental evidence on the sharp orientation tuning of V1 neurons (see also \cite{bresscow03} and \cite{bresscow02}).
 The selection of a direction $X_3$ defines at each point a subspace of dimension 2 of the Euclidean tangent space to $\mathbb{R}^2\times S^1$, generated by
 \begin{equation*} X_1 = \cos\theta \partial_x + \sin\theta \partial_y,\quad X_2 = \partial_\theta. \end{equation*}
 The planes generated by $X_1$ and $X_2$ at each point are called \emph{horizontal planes} and determine a bracket-generating distribution on $\mathbb{R}^2\times S^1$, thus defining a sub-Riemannian structure on it. In particular one has
 \begin{equation}
  [X_1,X_2] = - X_3,
 \end{equation}
 so that the Lie algebra generated by $X_1$ and $X_2$ is the whole three-dimensional tangent space at every point. This property implies, thanks to the Chow Theorem (see e.g. \cite{montgomery}), that any two points $(x_0,y_0,\theta_0)$ and $(x,y,\theta)$ of $\mathbb{R}^2\times S^1$ can be connected through an integral curve of $X_1 + k X_2$, where $k$ varies in $\mathbb{R}$. Such curves are called \emph{horizontal curves}.\\
 \begin{figure}[!hbtp]
\centering
 \includegraphics[width=0.5\textwidth]{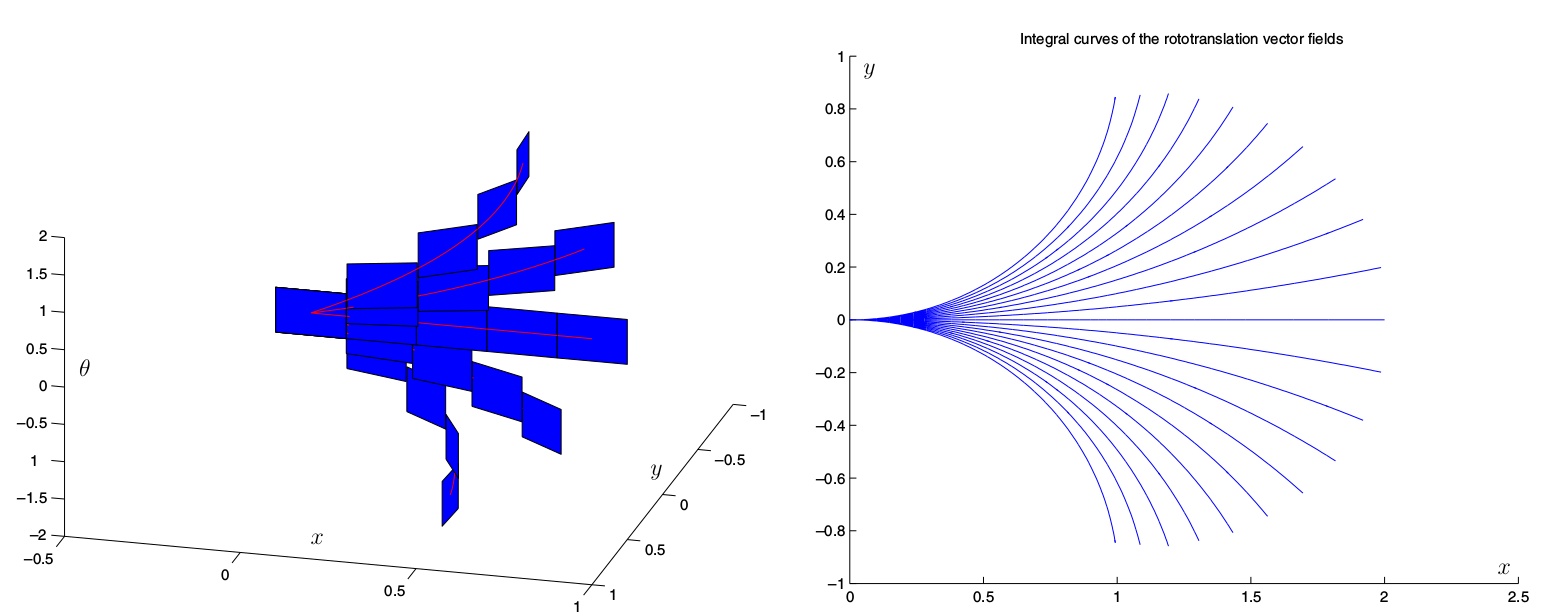}
 \caption{Integral curves of $X_1 + k X_2$ for varying values of $k$. Left: the curves visualized in $\mathbb{R}^2 \times S^1$, with contact planes displayed. Right: their projection onto the image plane. Source: \cite{cs06}.}\label{intcurves}
\end{figure}
 Just as the anatomical mechanism of horizontal connections is proposed to implement the perceptual association fields, the connectivity defined by this sub-Riemannian structure allows to recover such a pattern: the family of horizontal curves starting at a point $(x_0,y_0,\theta_0)$ can be interpreted as the local association field spreading from $(x_0,y_0,\theta_0)$, see Figure \ref{intcurves}.\\ 
 The propagation of neural activity in this sub- Riemannian space can be modeled by means of a diffusion equation expressed through second order operators defined in terms of the vector fields $X_1$ and $X_2$. A solid reference for these topics is the book \cite{blu}.
 \begin{figure}[!hbtp]
\centering
 \includegraphics[width=0.3\textwidth]{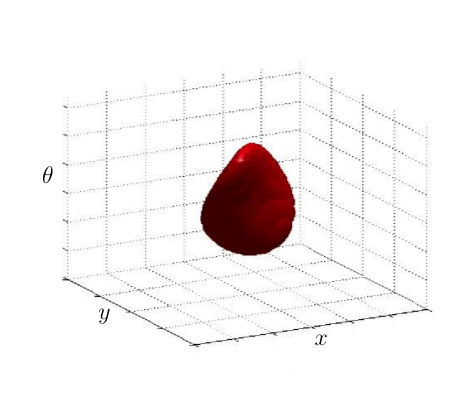}
 \caption{A level set of the fundamental solution of the sub-Riemannian Laplacian.}\label{lapl}
\end{figure}
 In particular, the heat equation associated to the sub-Riemannian Laplacian $\Delta u = X_1^2 u + X_2^2 u$ is considered. Figure \ref{lapl} displays a level set of the fundamental solution to this operator.\\
 
 The local invariance of the architecture of V1 w.r.t. $SE(2)$ has been exploited in a number of other works. In \cite{uncertainty}, cortical orientation maps are modeled as minimizers of an uncertainty principle stated in this setting. An extension of \cite{cs06} in the context of spatio-temporal analysis of visual stimuli was developed in \cite{bccs}. On the other hand, a semi-discrete variant of \cite{cs06} was proposed in \cite{boscain14}, where only a finite number of angles is taken into consideration. Moreover, a further study on the relation between association field curves and sub-Riemannian geodesics in $\mathbb{R}^2 \times S^1$ was carried out in \cite{duitsboscain}.\\
 In \cite{edge-stat}, a Fokker-Planck equation was considered instead of a sub-Riemannian heat equation, resulting from a study of statistical kernels of edge co-occurrence in natural images (see also \cite{mumford} and \cite{perceptual}). This corresponds to assuming that the propagation has a deterministic component along $X_1$ and a stochastic component in direction $X_2$.\\
 Linear as well as non-linear diffusion equations in $\mathbb{R}^2\times S^1$ were also employed in \cite{duitsI,duitsII} for medical image processing; here, the lifting process of the image is based on a so-called orientation score, essentially providing a local orientation representation of an image through a generalized unitary wavelet transform \cite{duits-score}. In \cite{abbfav}, such a lifting was complemented with the definition of a sub- Riemannian structure in a 5-dimensional space obtained by introducing the features of curvature and intensity, in order to overcome limitations given by crossing and bifurcations in a blood vessel tracking task. In \cite{duitsfuehr}, the approach of evolution equations was extended to the range of Gabor transforms, again with applications to the processing of images.\\
 
 Another model of V1 horizontal connections based on differential geometry has been proposed in \cite{benzuck}, where a crucial role is again given to curvatures, providing a directional rate of change of orientation. In this work, the relationship between nearby tangents is analyzed for curves as well as textures. An analysis of other possible relationships between the lateral connectivity and visual function beyond contour integration had been carried out in \cite{benshahar} as well.\\
 Moreover, a number of models have been proposed where perceptual grouping and texture segmentation tasks were addressed through a combination of long-range interactions and recurrent processing, involving the so-called bipole model \cite{gromin,neumin,hansneu}.\\
  
 In the next sections, we will propose a new model of the functional architecture of V1 that does not need the differential structures of most previous models. It is based on a metric structure on a suitable \emph{feature space} indexing the RPs of simple cells. Such a geometry is induced by the RPs themselves. Interestingly, when the family of RPs is represented by a bank of Gabor filters indexed by $\mathbb{R}^2 \times S^1$, the induced distance is locally equivalent \cite{metric} to a Riemannian approximation to the sub-Riemannian structure described in \cite{cs06}. Still, our model can be applied to a wide range of filter banks, with no need for any group structure, and is able to recover the perceptual rules of association fields even when applied to a randomly ordered family of learned orientation-selective filters.

\section{A functional architecture generated by RPs}\label{model}

\subsection{The generating kernel}\label{generator}

 We fix a bank of real- or complex-valued linear filters on the plane, modeling a set of RPs of simple cells of V1. Specifically, we consider a family $\{\psi_p\}_{p\:\in\:\mathcal{G}} \subseteq L^2(\mathbb{R}^2)$, indexed by a \emph{feature space} $\mathcal{G}$. These can be taken to be Gabor filters as well as any other bank of filters, for instance obtained through automatic learning procedures or to fit the shape of a set of experimentally measured RPs. As we will show later, endstopped simple cells may be included as well by adding a suitable sub-family of filters. Typically, the feature space has the form $\mathcal{G}=\mathbb{R}^2\times\mathcal{F}$, where $(x,y) \in \mathbb{R}^2$ represents the retinal position at which the filter is centered and $f \in \mathcal{F}$ encodes the other features detected by the filters. Let us define, for every $p,p_0 \in \mathcal{G}$, the \emph{generating kernel}
 \begin{align*}
 K(p,p_0) := Re\langle\psi_p,\psi_{p_0}\rangle_{L^2} = Re \left(\int_{\mathbb{R}^2} \psi_p(x,y) \: \overline{\psi_{p_0}(x,y)} \: dx \: dy \right).
 \end{align*}
 \begin{equation}\label{defK}\end{equation}
 This kernel is the reproducing kernel \cite{rk} on the closure of the linear span of $\{\psi_p\}_p$ in $L^2(\mathbb{R}^2)$ (see also \cite{duitsRK}); it extends to a general family of filters the familiar idea of the reproducing kernel associated to a family of wavelets \cite{antoine,deng}. In the special case where $\mathcal{G}$ has a group structure and the filters are given by $\psi_p = \mathcal{U}_p\psi$ for a fixed \emph{mother filter} $\psi$ through a unitary group representation $\mathcal{U}$, one has 
 \begin{equation}\label{groupcase}K(p,q)= K(p^{-1}q, \: \mathbf{1}).\end{equation}
 This is indeed the case for the example of a bank of Gabor filters that will be examined later.\\
 
 The kernel $K$ expresses the strength of \emph{correlation} between the RPs of two cells, represented by the points $p$ and $p_0$ of the feature space. Now, it is not restrictive to assume that the filters are normalized to have squared $L^2$-norm equal to some $\eta>0$. Under this assumption, we define the \emph{kernel distance}
 \begin{equation}\label{dK}
   d^2(p, p_0)  = 2\big(\eta - K(p,p_0)\big).
\end{equation}
 In fact, $d^2$ coincides with the squared $L^2$ distance between the functions $\psi_{p}$ and $\psi_{p_0}$. Indeed, $\| \psi_p - \psi_{p_0} \|_{L^2}^2$ can be rewritten as
\begin{align*}
   \|\psi_{p}\|_{L^2}^2 + \|\psi_{p_0}\|_{L^2}^2 - 2Re\langle\psi_p,\psi_{p_0}\rangle_{L^2}.
\end{align*}

 The generating kernel $K$ and the distance $d$ have a local sense. Indeed, the filters typically have a very localized support and it does not make much sense to compute $L^2$ scalar products between them if they are not sufficiently close together in the space domain. Moreover, the shape of the specific filters taken into consideration may bring one to introduce some restriction on which filters can directly be influenced by one another. This corresponds to defining around each point $p_0 \in \mathcal{G}$ a local distance function defined on a \emph{patch}, which we shall call $\mathcal{P}(p_0)$, and to gluing all these local distances together to obtain a global distance function
 \begin{align}\label{distinf}
 \tilde{d}(p,p_0) = \inf \Big\{ \:\sum_{j=1}^{N} d(q_{j-1},q_j) \: : \: N \in \mathbb{N},\: q_0=p_0,\:q_N=p,\:q_j \in \mathcal{P}(q_{j-1}) \:\forall j \Big\} .
 \end{align}
 Note that the length $N+1$ of the sequence $\{q_j\}$ is let vary: the $\inf$ is taken over all $N \in \mathbb{N}$ and over all sequences compatible with the patches.
 Under a non-degeneracy condition on the definition of such patches, this is still well-defined as a distance function (see \cite{metric} for more details).\\
 
 Although we restrict our analysis to the pairwise relations between RPs, some studies have emphasized the significance of higher order components. In \cite{lawzuck}, third order statistics of natural scenes are examined: in terms of reciprocal influences between cells, this would mean characterizing not only the probability of each cell to be activated given the activation of a ``central'' neuron, but also the probability of \emph{co-activation} of two neurons given the activation of the central one. In geometrical terms, such a model allows for curvature dependencies to emerge. See also \cite{benzuck} and \cite{abbfav}. As a matter of fact, we will show that some information on the curvature tuning of some classes of neurons, namely endstopped simple cells, emerges in our model as well.
 
 \subsubsection{Example: Gabor filters.}
 
  Let us now consider, as a first basic example, the case of a bank of Gabor filters of fixed scale $\sigma$. Its general element reads: 
 \begin{equation} \label{gaborfilt} 
 \psi_{x,y,\theta} (u,v) = \exp\left(\frac{2\pi i X}{\lambda}\right) \exp\left(-\frac{X^2+Y^2}{2\sigma^2}\right),\end{equation}
 with \begin{equation*} \begin{cases} X = (u-x)\cos\theta + (v-y)\sin\theta \\ Y = -(u-x)\sin\theta + (v-y)\cos\theta.\end{cases} \end{equation*}
 Denote $p = (x,y,\theta)$. A direct calculation yields:
\begin{align*} K(p,0) =  \sigma^2\pi \: \exp\left( -\frac{x^2}{4\sigma^2} -\frac{y^2}{4\sigma^2} - \frac{2\sigma^2\pi^2 (1 - \cos\theta)}{\lambda^2} \right) \cdotp  \cos\left( \pi \: \frac{x(1 + \cos\theta) + y\sin\theta}{\lambda} \right). \end{align*}
 The distance induced by $K$ is given by:
 \begin{align}\label{distgabor}
d^2 (p,0) = 2\sigma^2\pi - 2 \: K(p,0).
\end{align}
 Recall that the feature space in this case is $\mathcal{G} = \mathbb{R}^2\times S^1 = SE(2)$: as pointed out before, the explicit expression for $K(p,p_0)$ for $p_0 \neq 0$ can then be obtained from the one for $K(p,0)$ through the group operation. We have:
\begin{align*}K \big( (x,y,\theta), \; (x_0,y_0,\theta_0) \big) = K\big( (R_{\theta_0}T_{(x_0,y_0)}(x,y),\: \theta-\theta_0), \; (0,0,0) \big). \end{align*}
 In \cite{metric} we show that, in the Gabor case, $d$ is estimated by a Riemannian distance, and that such a metric turns out to be a Riemannian approximation to the sub-Riemannian structure defined in \cite{cs06} on $\mathbb{R}^2\times S^1$. Specifically, such convergence takes place when $\sigma^2$ is proportional to the wavelength $\lambda$ and $\lambda \rightarrow 0$. This means that the number of oscillations of the filters under the  Gaussian bell goes to infinity; on the other hand, their (squared) norm $\eta=\sigma^2 \pi$ shrinks.

\subsection{Propagation by repeated applications of $K$}\label{propagation}

A cortex model must be able to explain not only the action of RPs but also contextual influences and the emergence of association fields. Several existing differential models (reviewed in Section \ref{previous}) define on the feature space a Lie group structure, on which the horizontal connectivity is described by means of a diffusion equation expressed through second order operators associated to an invariant sub-Riemannian metric. By contrast, in the present model the local geometry of the cortical space is determined by the receptive profiles of the neurons themselves, through a generating kernel which expresses the reciprocal influences between simple cells with overlapping RFs. The long-range connectivity is then described by iterating the action of this kernel.\\

 It is important to remark that the horizontal connections take place in V1 layer 2/3, which is composed of phase-invariant complex cells. Note that
\begin{align}\label{evenodd} K(p, p_0) = Re\langle \psi_p, \psi_{p_0} \rangle_{L^2} = \langle Re\psi_p, Re\psi_{p_0} \rangle_{L^2} + \langle Im\psi_p, Im\psi_{p_0} \rangle_{L^2},\end{align}
 which for $p=p_0$ reads
\[ K(p,p) = \|Re\psi_p\|_{L^2}^2 + \|Im\psi_p\|_{L^2}^2.\]
 This recalls the classical energy model for V1 complex cells. In the case of Gabor filters, this is indeed the square sum of a quadrature pair of filters, since the real and imaginary parts of $\psi_{x,y,\theta,\sigma}$ show the same orientation and scale but are shifted in phase of $90^\circ$.\\
 Now, how does the activity spread to reach neurons with widely separated RFs? First, consider the operator
 \begin{equation}\label{step1}
  S[K](p,q) :=  \frac{K^{(1)}(p,q)}{\int_\mathcal{G} K^{(1)}(p,q')d\mu(q')}
 \end{equation}
 with
 \[K^{(1)}(p,q) = \frac{h(K(p,q))}{\int_\mathcal{G} h(K(p',q))d\mu(p') \:\cdot \: \int_\mathcal{G} h(K(p,q'))d\mu(q')},\]
 where $h$ is a nonlinear function. Note that the kernel is first normalized w.r.t. both its variables, followed by a further \emph{anisotropic} normalization, as done in \cite{coiflaf}, where a kernel expressing the local geometry of a given dataset is defined. Then consider the new kernel computed around a point $p_0$, i.e.
 \begin{equation}\label{laf}K_1^{p_0}(p) := S[K](p,p_0).\end{equation}
 The simplest choice for the function $h$ is $h(z):=\max(z-\tau, \:0)$. Setting $\tau=0$ implies discarding the inhibitory effects (represented by negative values), while $\tau>0$ would mean only keeping the strongest excitatory responses, thus essentially implementing the principle of non-maxima suppression at every diffusion step (as in \cite{bresscow03}, \cite{bresscow02} or \cite{cs06}). Another common choice of $h$ is the \emph{logistic sigmoidal function}
 \[h(z) = \frac{1}{1+e^{-z}} \: .\]
 We now iteratively define, for each $n>0$,
 \begin{equation}\label{iteration} K_{n}^{p_0}(p) := \int_\mathcal{G} S[K](p,q) \: K_{n-1}^{p_0}(q) d\mu(q).\end{equation}
 Note that $K_n^{p_0}$ integrates to 1 for every $n$. Alternatively, one may apply the same iteration rule to the original kernel, i.e. take $K_1^{p_0}(p) := K(p,p_0)$. This would mean keeping the initial inhibitory information. In both cases, one has: 
 \begin{equation}\label{normalization} |K_n^{p_0}(p)| \leq \|K_{n-1}^{p_0}\|_\infty \underbrace{\int_\mathcal{G} S\left[K\right](p,q) d\mu(q)}_{=1} \quad \forall \; n.\end{equation}
 Also note that, in the presence of a group structure, (\ref{groupcase}) holds and the integral in (\ref{iteration}) becomes a group convolution.\\
 Here, $\mu$ is the spherical Hausdorff measure \cite{hausd} associated to the distance $d$ on $\mathcal{G}$. This choice is motivated by the fact that, in the general case, the feature space $\mathcal{G}$ is just a set of parameters with no geometric structure whatsoever, apart from the one given by the distance: since we are integrating a function of $d$, it is natural to do so with respect to a measure induced by the distance itself, which ``compensates'' for its in-homogeneities. In a simpler case, suppose that $\mathcal{G}$ is a subset of some $\mathbb{R}^n$, which is in fact equipped with a measure: the actual space might be much lower-dimensional, and integrating with respect to the $n$-dimensional Euclidean measure would make no sense in this case. On the contrary, the spherical Hausdorff measure captures the ``true'' dimension of the space.\\
 
 Note that the spatial extent of the kernel widens at each step, and this width depends on the size of the RPs, i.e. on the diameter of their RFs (RPs are local objects, and they are always modeled through compactly supported or rapidly decaying functions).\\
 There is neurophysiological evidence \cite{angelucci02} that the extension of horizontal connections departing from a V1 cell matches the size of its so-called \emph{low-contrast summation field}; this is identified as the area measured by presenting low-contrast bars or gratings of increasing sizes at the RF center, and keeping the size of peak response. The low-contrast summation field has in turn been shown \cite{kapadia,angelucci02} to be on average 2.2-fold greater than the ``classical'' RF. Therefore, since a relationship between the RF size and the spatial extent of horizontal connections seems to hold from a biological point of view, we believe that the optimal number of iterations may not depend on the RP size. As for the choice of this number, a first stopping criterion is this average ratio between the size of low-contrast summation fields and classical RFs. Another possible quantitative framework for fitting the kernel to the neurophysiological data was proposed in \cite{favali16b} to fit a Fokker-Planck kernel to the data from \cite{bosking} and \cite{angelucci02}: the kernel was evaluated on a pinwheel map and compared with the measured distribution of the tracer by comparing their densities onto the rectangles of a grid (whose sampling size was chosen to match the pinwheel scale).\\

 Note that this idea of repeated convolutions can be applied to model the response of the layers of V1 to a visual stimulus $I$ applied to the retina, taking into account the contextual influences. In this case, the activation computed by lifting the image to the cortex is a function of the cortical coordinates:
 \[ I_0(p) := h\left(\int I(x,y)\psi_p(x,y) dxdy\right). \]
 The action of the kernel on this activation at the $n$-th step may then be expressed by
 \begin{equation}\label{imevol} I_{n}(p) := \int_\mathcal{G} S[K](p,q) \: I_{n-1}(q) d\mu(q). \\\\\end{equation}

 In our work \cite{metric} we constructed a diffusion process on the feature space, equipped with the measure $\mu$, to describe the propagation along the connectivity. Subject to some requirements \cite{sturm95b,sturm98,coiflaf}, this process can be approximated by updating the initial datum through repeated integrations against a proper normalization of $K$: indeed, under certain conditions such a kernel provides an estimate for the heat kernel associated to the diffusion process \cite{coiflaf}. See \cite{metric} for further details. This remark highlights a link between our kernel-based model and the diffusion-based differential models cited before, showing that these two approaches can actually coexist.

\section{Numerical scheme and results}\label{results}

 In this section we show, through some numerical simulations, the geometrical properties of the connectivity kernels obtained by applying our model to some different banks of filters. The first distinction to be made in the treatment of these examples from the numerical point of view is that between the \emph{continuous} case and the \emph{discrete} one. This can refer either to the feature space $\mathcal{G}$ or to the spatial domain of the single filters. For instance, in the case of a Gabor system, both the feature space $\mathcal{G} = \mathbb{R}^2 \times S^1$ and the spatial variable $(u,v) \in \mathbb{R}^2$ need to be sampled. On the contrary, when dealing with a bank of learned filters, the family of indices $\mathcal{G}$ and the filter domains are finite sets. Of course, one could also consider ``hybrid'' cases, e.g. a bank of filters defined on $\mathbb{R}^2$ but indexed by a discrete set.\\ 
 
 For what concerns the computation of the initial kernel, the numerical setup is very basic. The filters are either defined on $\mathbb{R}^2$ or numerically known onto a set of the type $\{(\frac{i}{n},\frac{j}{n})\}_{i,j=-n,\ldots,n} \subseteq \mathbb{R}^2$; computing the generating kernel $K(p,q)$ simply involves taking the integral of the function $\psi_p \cdot \overline{\psi_q} \in L^1(\mathbb{R}^2)$, which is either compactly supported or exponentially decaying. Thus, the filter domains can safely be sampled (when needed) through a bounded square grid 
 \[\{(x_i,y_j)\}_{i,j=0,\ldots,N} \subseteq [-W,W]\times[-H,H],\]
 for some $N,W,H>0$, with 
 \[x_{i+1}-x_i = y_{j+1}-y_j =: \delta > 0 \quad \forall i,j.\]
 The kernel $K(p,q)$ is then standardly calculated as
 \[\frac{1}{\delta^2} \sum_{i,j} \psi_p(x_i,y_j) \overline{\psi_q(x_i,y_j)}.\]
 In the Gabor case, as shown in Section \ref{generator}, this integral can even be computed analytically, thus avoiding this approximation.\\
 As for the propagation, first note that, in the case of a finite feature space $\mathcal{G}=\{p_1,\ldots ,p_M\}$, the Hausdorff measure reduces to the counting measure: integrals over $\mathcal{G}$ are finite unweighted sums over its elements. Although seemingly trivial, this indeed captures the distribution of the features throughout the data thanks to the inhomogeneity of the distance $d$. Loosely speaking, given a certain set of filters, some features may be ``more densely'' represented than others w.r.t. $d$. As a basic example, suppose that the filter $\psi_{p_1}$ is highly correlated with 8 other filters, i.e. for a fixed $\varepsilon$
 \[ \big|\big\{i \in \{2,\ldots ,M\} \: : \: d(p_i,p_1)<\varepsilon\big\}\big| = 8.\]
 On the other hand, suppose that $\psi_{p_2}$ has only 2 filters highly correlated with it. As a consequence, $\mu(B_\varepsilon^d(p_1)) = 9$ and $\mu(B_\varepsilon^d(p_2)) = 3$, i.e. the balls of $d$ of the \emph{same} radius $\varepsilon$ centered at $p_1$ and $p_2$ have different measures. For a non-finite metric space, the spherical Hausdorff measure extends this concept, although in general it is not explicitly computed; however, the distance $d$ in the Gabor case turns out to be estimated by a Riemannian distance \cite{metric}, and the spherical Hausdorff measure on a Riemannian manifold coincides up to a constant factor with the Riemannian volume form \cite{federer}. We will go into more detail on the discretization of the Gabor feature space in Section 4.1.\\
 We stress that this work focuses on showing that our construction allows to recover a connectivity pattern starting from non-structured banks of filters. A great interest is thus addressed to numerically known (therefore finite) sets of filters, obtained for instance through learning algorithms; one may even apply it to neurophysiologically measured RPs. We refer to \cite{metric} for some theoretical issues and more examples related to the spherical Hausdorff measure in that context.\\
 A last point to be mentioned is the role of the normalization performed at each kernel iteration. Namely, as shown in Eq.(\ref{normalization}), $\|K_n^{p_0}\|_\infty$ is non-increasing w.r.t. $n$, which prevents the values of $K_n^{p_0}$ from exploding as $n \rightarrow \infty$.\\

 In the following, we first resume the classical example of Gabor filters presented in Section \ref{generator} and display the association fields generated by the connectivity in this case. We then show the curvature selectivity of the association fields emerging from endstopped simple cells. We also take into account the spatiotemporal behavior of the RPs by considering a family of three-dimensional Gabor filters including a time parameter. Finally, we show that it is still possible to recover a coherent ``bow-tie'' pattern starting from a family of RPs obtained through a learning algorithm, whose parameterization carries no a priori geometric information.

  \subsection{Gabor filters}\label{gabor}

 We now go back to the bank of Gabor filters $\{\psi_{x,y,\theta}\}$ defined by Eq. (\ref{gaborfilt}). Since the generating kernel in this case is known analytically, no discretization of the filter domain $\mathbb{R}^2$ is needed for its computation. As for the feature space $\mathbb{R}^2\times S^1$, in the following we consider a sampling of the type
 \[ \{-\overline{x},\ldots,\overline{x}\}\times\{-\overline{y},\ldots,\overline{y}\}\times\{-\overline{\theta},\ldots,\overline{\theta}\} \]
 where $\overline{x}=\overline{y}=1$ with discretization step .01 and and $\overline{\theta}=1.5$ with step .015 for the visualizations of the generating kernel; $\overline{x}=1.5$ and $\overline{y}=3$ with step .1, and $\overline{\theta}=1.5$ with step .15 for the propagation.\\
 
  \begin{figure}[!htbp]
 \centering
 \includegraphics[width=.4\textwidth]{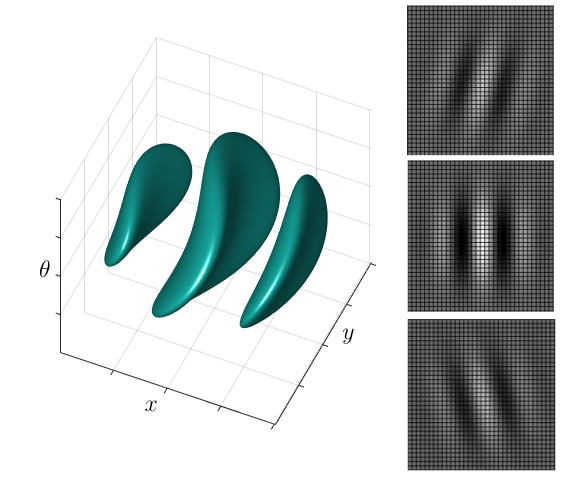}
 \caption{On the left, a level set of $K\big((x,y,\theta), (0,0,0)\big)$. On the right, three slices for $\theta = -\frac{\pi}{6}, 0, \frac{\pi}{6}$.}\label{levelset}
\end{figure}
 In the case of Gabor filters, the level sets of $K$ (or equivalently, of the distance $d$) around a fixed point $p_0=(x_0,y_0,\theta_0)$ of the feature space are in general not connected (see Figure \ref{levelset}). This is due to the oscillations of the periodic factor.
 \begin{figure}[!htbp]
 \centering
 \includegraphics[width=.4\textwidth]{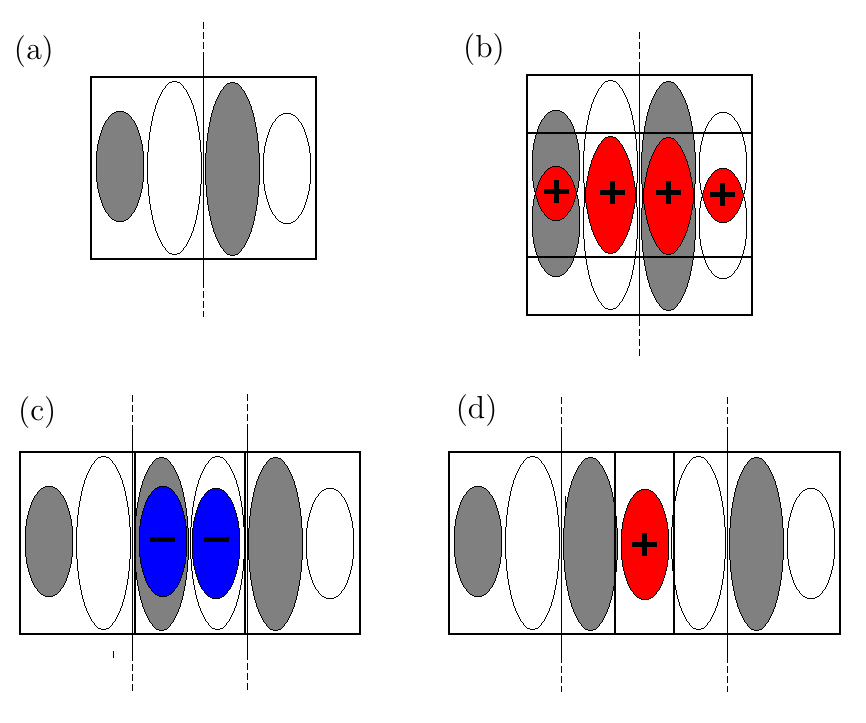}
 \caption{A schematic representation of the interactions between odd RPs sharing the same orientation, depending on their reciprocal position. (a) An example odd filter, with white ON areas and gray OFF areas. The axis of its preferred orientation is displayed. The next pictures show different reciprocal configurations of two such filters. Overlapping areas where the sign of the product is positive are displayed in red and marked by a plus symbol, while zones where the product is negative are blue and marked by a minus symbol. (b) Collinear filters bring to high values of $K$ along their common axis. (c) Parallel filters with overlapping regions of opposite signs yield negative values of $K$. (d) Parallel filters with overlapping regions of the same sign yield positive values of $K$.}\label{overlap}
\end{figure}
 \begin{figure}[!htbp]
 \centering
 \includegraphics[width=0.5\textwidth]{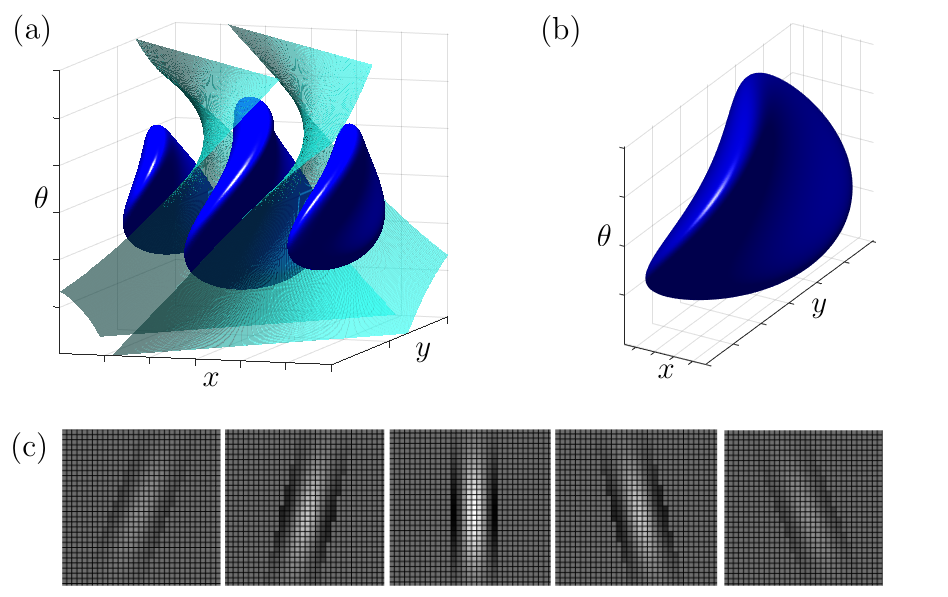}
 \caption{(a) In dark blue, a level set of $K\big((x,y,\theta), (0,0,0)\big)$. The patch $\mathcal{P}(0,0,0)$ is the area between the two surfaces displayed in light blue. (b) The same level set, after truncating. Here, $\lambda=1$. (c) Horizontal slices of the truncated kernel for $\theta=-0.75, -0.45,0,0.45,0.75$. The white regions represent positive values, the black ones represent negative values. Note that the kernel has been truncated at its minimum.}\label{patch}
 \end{figure}
 In terms of the interactions between RPs, the central lobe of the level sets corresponds to the \emph{co-axial} effect, i.e. the one along the axis of the preferred orientation of the starting RP $\psi_{p_0}$: along this axis, the strongest interactions are the ones with filters collinear with $\psi_{p_0}$ (see Figure \ref{overlap}(b)), and the value of $K$ decreases as the orientation of the filters varies from $\theta_0$. On the other hand, the smaller lobes developing along the orthogonal axis result from the ``periodic'' correlation of the central RP with filters parallel to $\psi_{p_0}$ (i.e. filters with orientation $\theta_0$ but centered on points of the orthogonal axis): the kernel takes positive values when areas of the two RPs with the same sign overlap (as in Figure \ref{overlap}(d)), and negative values when areas of opposite signs overlap (as in Figure \ref{overlap}(c)). In a first stage we thought it best to restrict ourselves to the effect along the axis of the preferred orientation, i.e. to consider only the central connected component of the level sets. This only takes a straightforward analytical operation in the Gabor case, and it makes it easier to compare our model with the ones obtained through real-valued diffusion equations. We shall therefore clip the kernel as anticipated in the general description above, according to local patches defined around each point $p_0=(x_0,y_0,\theta_0)$ of $\mathbb{R}^2\times S^1$:
 %\begin{small}
\[ \mathcal{P}(p_0) := \{ (x,y,\theta) \: : \: |a(1 + \cos\delta) + b\sin\delta| < \lambda \}, \]
 %\end{small}
 with $(a,b,\delta) = \big(R_{\theta_0}T_{(x_0,y_0)}(x,y),\: \theta-\theta_0\big)$.
 This provides an example of adjustment that can be introduced on the original kernel in order to define some restriction on which filters are influenced by one another. At this point we can truncate the distance function by defining the kernel $K$ only on $\mathcal{P}(p_0)$ around each point $p_0$ (see Figure \ref{patch}), and then glue the resulting distance functions together as in Eq. (\ref{distinf}) to obtain a new global structure.\\
 Also note that, since the reciprocal influence of collinear filters is greater that the one between parallel filters, all the values of $K(\cdot, p_0)$ above a sufficiently high threshold are confined inside $\mathcal{P}(p_0)$: this means that small balls of the distance $d$ are not altered by the clipping. In other words, the local structure is preserved.\\
 \begin{figure}[!htbp]
 \centering
 \includegraphics[width=0.3\textwidth]{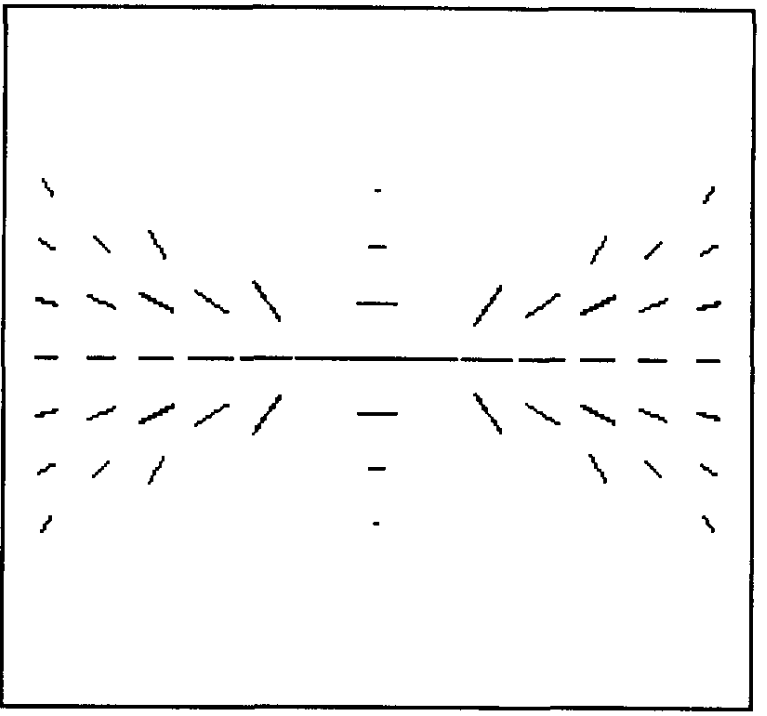}
 \caption{The \emph{connectivity pattern} found in the model of \cite{yenfinkel} includes both co-axial and trans-axial contributes. The length of the lines indicates connection strength.}\label{ladder}
 \end{figure}
 Nonetheless, we believe that the whole kernel is interesting for a further analysis, since its oscillatory behavior in the orthogonal direction seems to account for the \emph{trans-axial} ``ladder'' effect which has indeed been observed in both neurophysiological and psychophysical studies \cite{mitchcrick,field,yenfinkel} (see also Figure \ref{ladder}), thus naturally including in the model the Gestalt perceptual principle of parallelism.\\
 
 The generating kernel $K$ computed around a point $(x_0,y_0,\theta_0) \in \mathcal{G}$ is a three-dimensional function
  \[(x,y,\theta) \mapsto K\big((x,y,\theta), \: (x_0,y_0,\theta_0)\big).\]
  Since the feature space for this family is of the form \emph{position} $\times$ \emph{orientation}, it is possible to display a projection of it on the retinal plane.
 \begin{figure}[h]
 \centering
 \includegraphics[width=.5\textwidth]{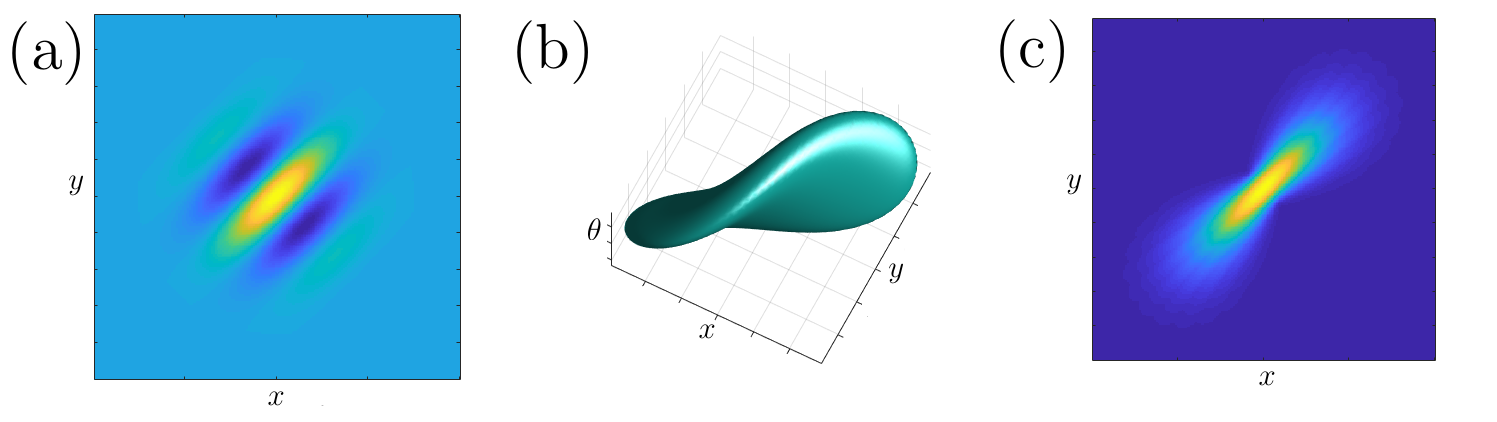}
 \caption{The behavior of $(x,y,\theta) \mapsto K\big((x,y,\theta), \: (0,0,\frac{\pi}{4})\big)$ in the case of Gabor filters. The kernel has been truncated as explained above. In particular, the images show: (a) the receptive profile $\psi_{p_0}$ corresponding to $p_0 = (0,0,\frac{\pi}{4})$; (b) a level set in $\mathbb{R}^2 \times S^1$; (c) the projection on the $(x,y)$ plane obtained taking the maximum in $\theta$. }\label{laf1}
 \end{figure}
 Specifically, by taking the maximum in the variable $\theta$ and projecting onto the $(x,y)$ plane, we obtain a 2D function concentrated around $(x_0,y_0)$, as displayed in Figure \ref{laf1}(c).\\
 
  We now display some stages of the iterative process described in Section \ref{propagation}, for the truncated kernel $K$.
\begin{figure*}[htbp!]
  \centering
  \includegraphics[width=0.7\textwidth]{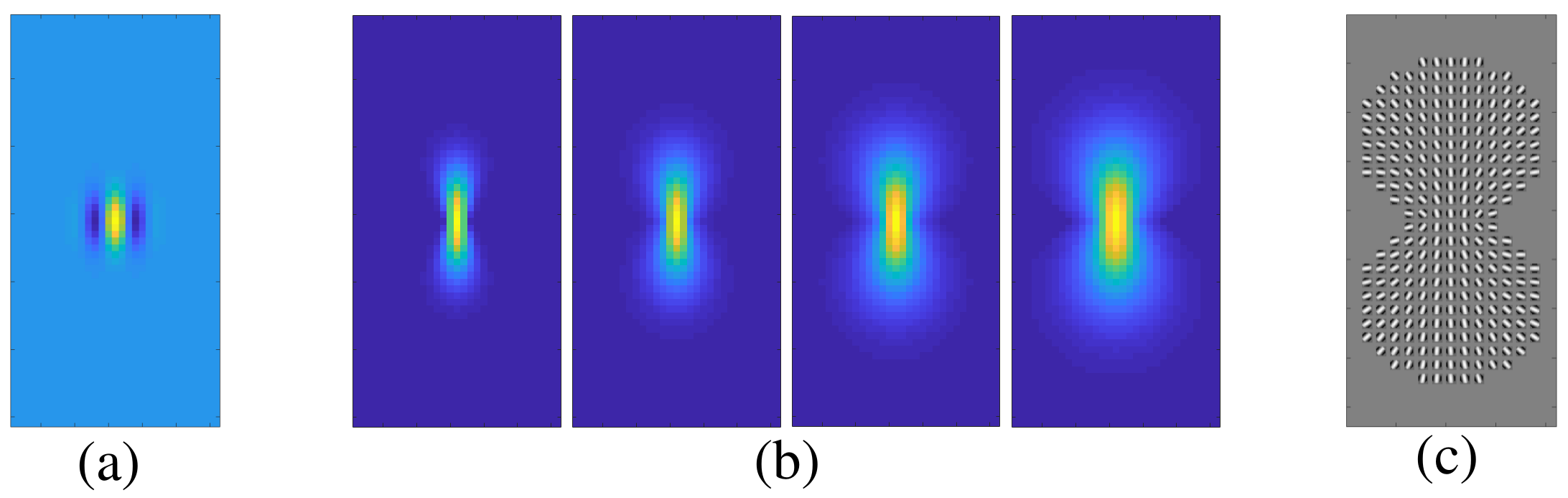}
  \caption{Development of a connectivity kernel around $(x_0,y_0,\theta_0)=(0,0,0)$, obtained through repeated integrations against the generator $K$. Four iterations were performed. The kernel was truncated as explained in Section \ref{generator}, and the operator $S$ defined in Section \ref{propagation} was applied. The resulting function $K^0_4(x,y,\theta)$ was then projected on the image plane by taking the maximum in the variable $\theta$. The pictures show: (a) the Gabor filter $\psi_{0,0,0}$ (real part); (b) the projection on the $(x,y)$-plane of the 3D function obtained at each of the three steps; (c) for each $(x,y)$, the value $\overline{\theta}=\arg\max_\theta K_4^0(x,y,\theta)$ represented through a Gabor filter with orientation $\overline{\theta}$.}
  \label{diffK}
\end{figure*}
 Figure \ref{diffK}(b) displays the projection onto the retinal plane of the kernel around $(0,0,0)$, at four steps of the iterative process.
 %The evolution has been implemented on $]-1,5,1.5[\times ]-3,3[ \times ]-1.5,1.5[ \subseteq \mathbb{R}^2 \times S^1$, discretized with step 0.1 in $x$ and $y$ and with step 0.15 in $\theta$.
 For $n = 1,2,3,4$ the function 
 \begin{equation}\label{Kn} (x,y,\theta) \mapsto K_n^{(0,0,0)}(x,y,\theta) \end{equation}
 was projected onto the $(x,y)$ plane by taking the maximum in the variable $\theta$. The projection of the connectivity kernel onto the retinal plane highlights the orientation selectivity of the propagation, with the points of maximum intensity concentrating along the axis of the preferred orientation of the central filter. As specified earlier, the horizontal connections of V1 are considered as a likely anatomical background to the perceptual rules synthesized by association fields. To develop this aspect, we added a further visualization displaying the orientation of the filter maximizing the intensity at each location (where the intensity exceeds a threshold). Specifically, Figure \ref{diffK}(c) shows at each such location $(x,y)$ a tiny Gabor filter whose orientation value $\overline{\theta}(x,y)$ is the one that maximizes the value of the final 3D function, i.e.
 $$ \overline{\theta}(x,y) := \arg\max_\theta K_4^{(0,0,0)}(x,y,\theta). $$
 
 Similar results have been obtained in several works concerning second-order statistics on the distribution of edges in natural images. In \cite{augzuck}, the statistical measurement of the organization of local edge elements into curves led to the introduction of an oriented filtering operation, expressed by a kernel on $\mathbb{R}^2 \times S^1$ (see Figure 15.3 in \cite{augzuck}); its orientation specificity and the negative lateral lobes recall some characteristics of our connectivity kernel, especially after clipping it at its minimum as explained above (cf. the slices in Figure \ref{patch}(c)). In \cite{geisler}, association field-like structures comparable to our Figure \ref{diffK}(c) emerge from both a statistical analysis of edge co-occurrence in natural images, and a Bayesian study analyzing the probability of two edge elements belonging to the same extended contour. A similar analysis is carried out in \cite{eldgold}, except that here contours are described as \emph{ordered} sequences of tangent elements; in this work a stronger role of the parallelism cue emerges, i.e. there is a marked correlation between parallel edge elements. See also \cite{kruger} and \cite{sigman} as further references on the emergence of the collinearity, cocircularity and parallelism cues from statistics on natural images.\\
 These studies are not to be regarded as a separate kind of analysis leading to results consistent with ours. Indeed, it is believed that the cortical processing is influenced by environmental statistics and adapts to best process those signals that are most likely to occur.

 \paragraph{Propagation on the pinwheel surface.}
 Consider now the subset $\{\tilde{\psi}_{x,y}\} = \{\psi_{x,y,\theta(x,y)}\}$ of the above family determined by a fixed orientation map $\theta(x,y)$.
 As pointed out in \cite{petitot}, orientation map-like structures can be reproduced as a superposition of plane waves with random phases. The generating kernel $\tilde{K}\big( (x,y),(x_0,y_0) \big)$ associated to this family is given by
 \[K\big( \: (x,y,\theta(x,y)), \; (x_0,y_0, \theta(x_0,y_0)) \: \big).\]
 \begin{figure}[!htbp]
  \centering
  \includegraphics[width=0.5\textwidth]{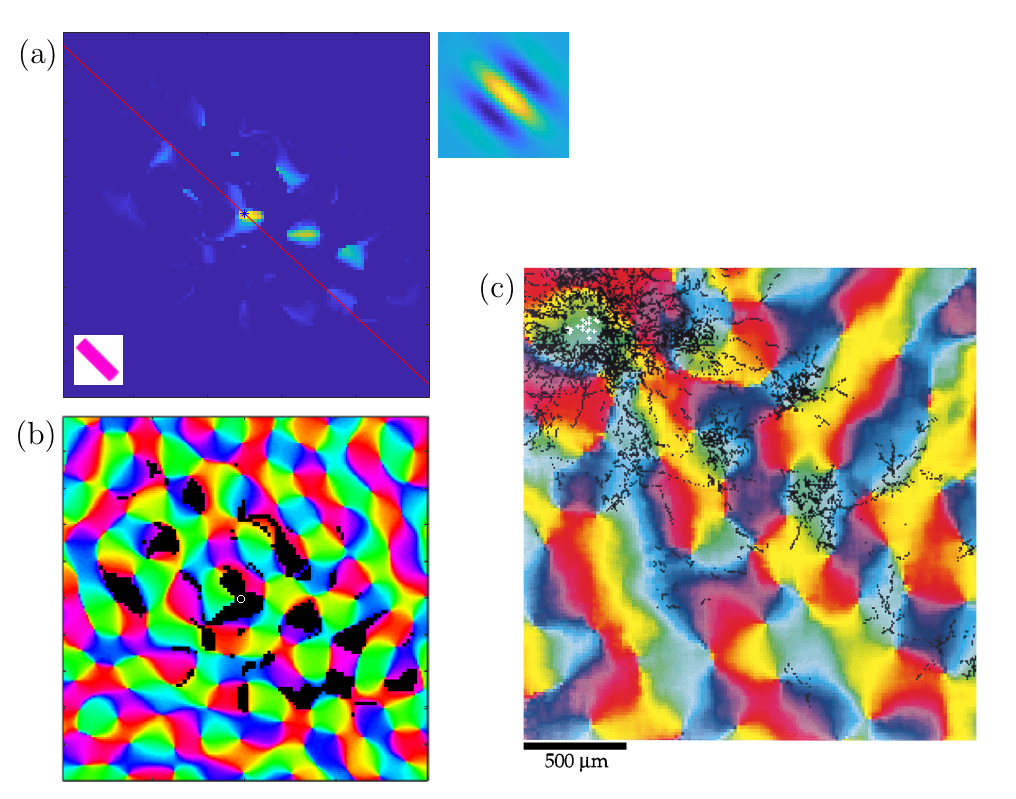}
  \caption{Left side: propagation onto a pinwheel map $\theta(x,y)$ generated through superposition of plane waves, as explained in \cite{petitot}. In particular, image (a) shows the function $\tilde{K}^{(0,0)}_6(x,y)$ resulting from six steps of evolution through $\tilde{K}$, with initial point $(0,0)$. In the lower left corner, a bar with orientation $\theta(0,0)$ is displayed, whose color corresponds to $\theta(0,0)$ in the orientation map. The starting filter $\tilde{\psi}_{0,0}$ (real part) is displayed as well (top right). Image (b) shows the values of $u$ exceeding a threshold, displayed in black onto the orientation map. The initial point $(0,0)$ is marked with a white circle. Right side: distribution of boutons formed after a biocytin injection into a site of tree shew V1, displayed over the topographic map of orientation preference. The cells that took up the biocytin are displayed in white. (Image (c), modified from \cite{bosking})}\label{pw}
  \end{figure}
  
  Starting from this 2D generator, we computed the connectivity kernel
  \begin{equation}\label{connk} (x,y) \mapsto \tilde{K}^{(0,0)}_n(x,y)\end{equation}
  onto the pinwheel surface $\theta$ parameterized by $(x,y)$, around its point $(0,0)$. We then displayed it (for $n=6$), as well as the orientation map itself, onto the $(x,y)$ plane, see Figure \ref{pw}(a). Note that, along with the orientation specificity already observed in the preceding experiment, a patchy nature of the propagation emerges as well, according to the regions of the map with orientation values close to $\theta(0,0)$. This behavior matches the one observed in \cite{bosking}, where small biocytin injections in a site of tree shew V1 were made to analyze the distribution of the resulting horizontal connections with respect to the orientation preference map. Figure \ref{pw}(c) shows the distribution of boutons obtained from one such injection. For the sake of comparison, we displayed onto the orientation map $\theta(x,y)$ the points (in black) where (\ref{connk}) exceeds a threshold (Figure \ref{pw}(b)). In both images, the starting point is displayed in white.\\
  In fact, in the experiment of \cite{bosking}, the recorded connections departing from the injection site originate not from a single cell, but rather from a few cells, the ones which took up the biocytin. However, all of these cells are confined in a localized area of V1: thanks to the continuity (outside pinwheels) of the orientation map, they happen to have very close preferred orientations. We believe that, if the byocitin injection could address a single neuron, the effect (although weaker) would be substantially the same in terms of orientation specificity. On the other hand, the multiple neurons situation may be reproduced in the context of our model by considering a set of RPs, close together in the feature space, and by combining the kernels generated by each of them (e.g. by summation).

  \subsection{Endstopped simple cells}\label{endstop}
  
  We will now focus on the connectivity pattern emerging from a particular class of neurons, namely \emph{endstopped simple} (ES) cells. Endstopped cells (also referred to as \emph{hypercomplex}) are orientation-selective neurons characterized by their sensitivity to bars of specific lengths: their RPs are characterized by inhibitory ``end zones'' suppressing the response to stimuli longer than the central excitatory area. Endstopped cells have been shown \cite{HW} to react to curved stimuli; a mathematical description of the relation between endstopping and curvature has been developed in \cite{dobbins87}. \\
  ES cells have \emph{simple} RPs equipped with end zones along the preferred orientation. They can be modeled by linear combination of two similarly positioned and oriented simple RPs of different sizes, as proposed in \cite{dobbins87}. Specifically, one such RP $\psi$ can be written as the (weighted) sum of the positive contribute of a smaller simple RP $\psi^S$ and the negative contribute of a larger simple RP $\psi^L$, as follows:
  \begin{equation}\label{endcomb}
   \psi := c_S \: \psi^S - c_L \: \psi^L,
  \end{equation}
  where $c_S > c_L >0$.
       \begin{figure}[htbp!]
 \centering
  \includegraphics[width=0.5\textwidth]{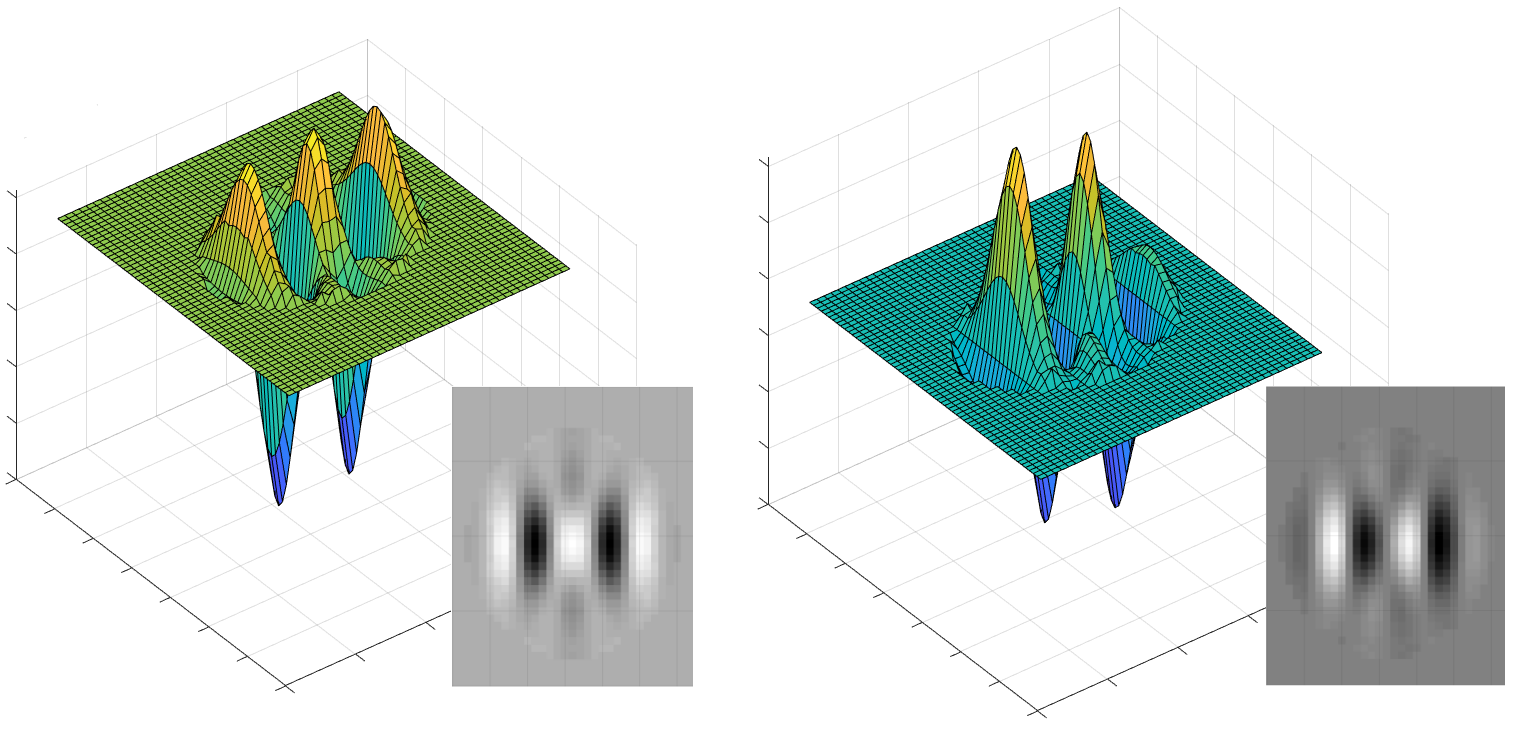}
  \caption{ES profiles obtained as in Eq. (\ref{endcomb}).}\label{endstoppedRPs}
 \end{figure}
 Figure \ref{endstoppedRPs} shows an example of even and odd RPs generated as in (\ref{endcomb}), where the simple RPs $\psi^S$ and $\psi^L$ are modeled by Gabor filters. To be more precise, the model in \cite{dobbins87} proposes to represent the response $R$ of $\psi$ to some contrast pattern as
 \begin{equation}\label{endresponse}
  R = h \big( c_S \: h(R_S) - c_L \: h(R_L) \big),
 \end{equation}
 where $R_S$ and $R_L$ are the responses of $\psi^S$ and $\psi^L$ to such pattern, and $h$ is the rectifying function $h(z) = \max(0, z)$. This increases the stimulus specificity compared to directly computing the response from (\ref{endcomb}).
     \begin{figure}[htbp!]
 \centering
  \includegraphics[width=0.6\textwidth]{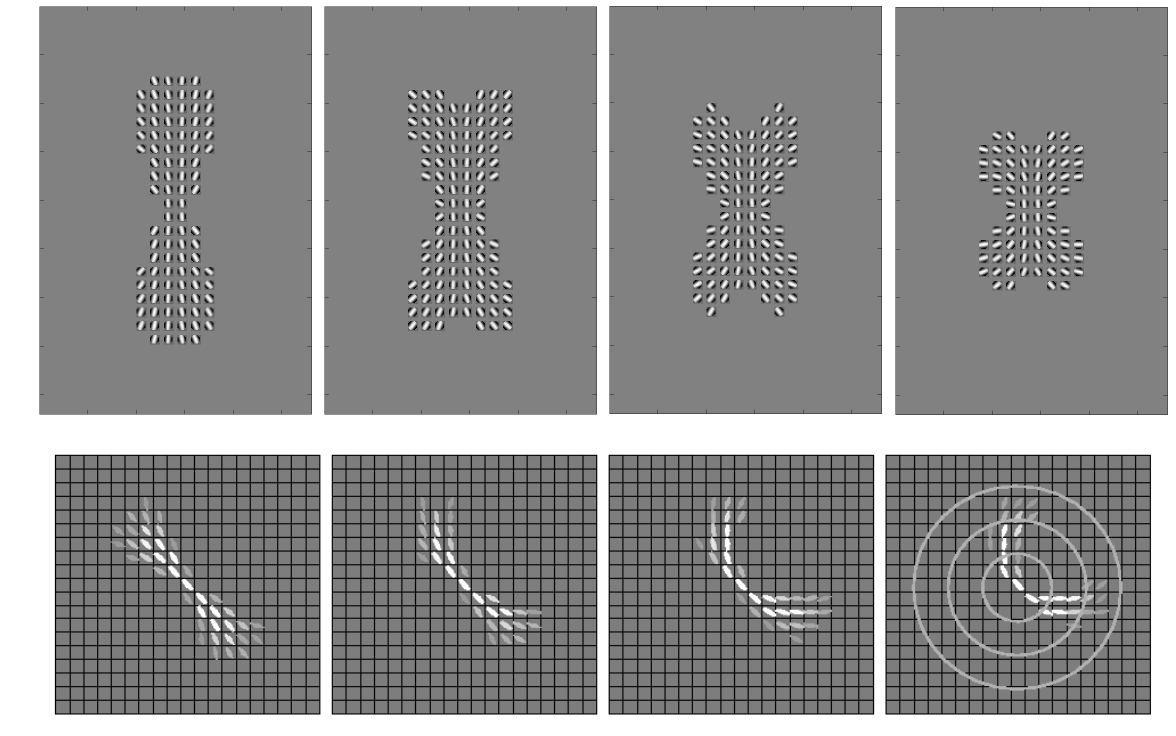}
  \caption{Top row: association fields emerging from the connectivity kernels associated to a simple cell (left) and to three ES profiles of varying length. Bottom row: the connection fields of curves obtained in \cite{benzuck} relative to curvature $\kappa=0$ (left) and three increasing positive values of $\kappa$ respectively.}\label{curvature}
 \end{figure}
 Nonetheless, the two models coincide when the responses of both $\psi^S$ and $\psi^L$ to a visual stimulus are nonnegative, i.e. for those stimuli in the ``preferred range'' of orientation and curvature of the ES cell.\\ 
 We computed the kernels associated to some families of ES profiles of different lengths and displayed the resulting association fields, obtained as before by selecting the orientation that maximizes the intensity at each location. This experiment shows that our model allows to recover the link between the length of ES cells and the curvature of the corresponding association fields. The top row of Figure \ref{curvature} shows the connectivity patterns obtained for a set of non-endstopped simple cells (left), corresponding to zero curvature, and for three families of ES cells of decreasing lengths, corresponding to increasing absolute values of curvature. In each of the cases displayed, we considered 21 RPs with orientation values equally spaced in $[-\pi,\pi]$. We then shifted them spatially to generate the whole filter bank. Our results are compared to the curvature-based connection fields obtained in \cite{benzuck} through differential tools, displayed in the bottom row of Figure \ref{curvature}. Since the sign of the curvature is taken into consideration in \cite{benzuck}, each of our association fields corresponds to a couple of theirs (with curvatures $\kappa$ and $-\kappa$).

 \subsection{Spatiotemporal Gabor filters}\label{spatiotemp}

 Let us consider a family of filters still of Gabor type, but taking into account the movement of the stimulus. Cocci et al. \cite{cocci} fitted the RPs of a set of V1 neurons showing velocity-selective behaviours with a three dimensional Gabor model. That is we have, in addition to the two spatial dimensions, a third temporal dimension in the domain of these filters, which in fact form a subset of $L^2(\mathbb{R}^3)$. 
 \begin{figure}[ht]
 \centering
  \includegraphics[width=0.5\textwidth]{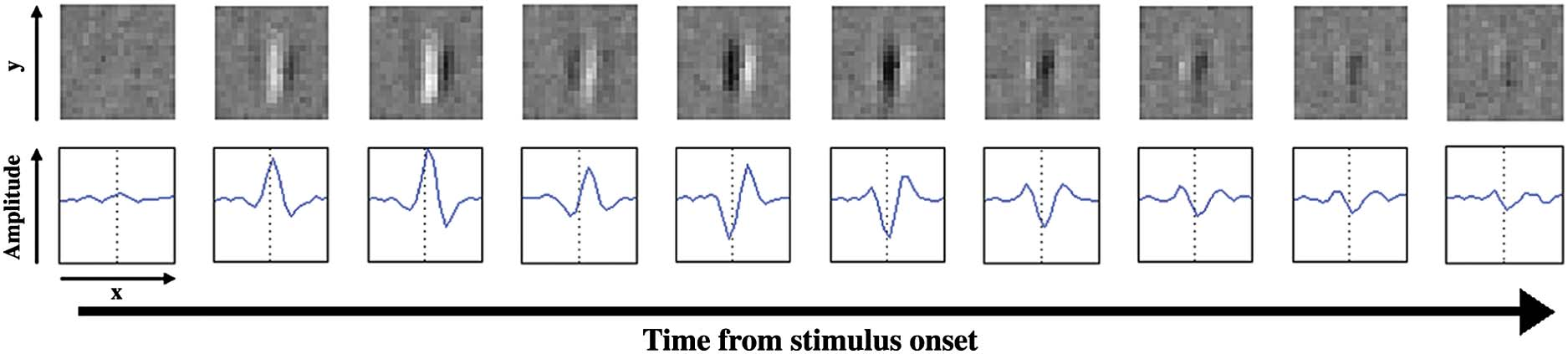}
  \caption{The time course of the recording of a simple cell's RP. Source: \cite{cocci}. }\label{spt}
 \end{figure}
 A convenient visualization of such filters is as a temporal sequence of spatial maps (see Figure \ref{spt}, taken from \cite{cocci}). In \cite{bccs}, the authors also develop a differential model of functional architecture based on such 3D Gabor family.\\
 Now, in order to mantain the notations as similar as possible to the one that we used for 2D Gabor filters, we shall write the general \emph{inseparable} element of the family of filters as 
 \begin{align}\label{gabortemp}\psi_{x,y,\theta,t,\alpha} (u,v,s) =  \exp\left(- 2\pi i \left( \frac{ X}{\lambda} + \alpha(s - t)\right)\right) \cdotp \exp\left(-\frac{X^2+Y^2}{2\sigma^2} - \frac{(s - t)^2}{2\beta^2}\right),\end{align}
 where \begin{equation*} \begin{cases} X = (u-x)\cos\theta + (v-y)\sin\theta \\ Y = -(u-x)\sin\theta + (v-y)\cos\theta.\end{cases} \end{equation*}
 The set of filters is indexed by the position $(x,y) \in \mathbb{R}^2$, the orientation $\theta \in S^1$, the time $t \in \mathbb{R}^+$ at which the response of the filter is maximum, and the velocity $\alpha \in \mathbb{R}$. Parameters $\lambda$ and $\beta$ are fixed.\\
 Let us compute the generating kernel $K$ in this setting. The feature space in this case is $\mathbb{R}^2\times S^1 \times \mathbb{R}\times \mathbb{R}$. We compute $K\big((x,y,\theta,t,\alpha),(x_0,y_0,\theta_0,t_0,\alpha_0)\big)$ with $t=t_0$: this means considering two cells whose activation peaks at the same time. In other words, we fix $t$ and we take $\mathcal{G} = \mathbb{R}^2\times S^1\times \mathbb{R}$ as our feature space. We restrict to this case in order to be able to better interpret the results in terms of the orientation and velocity parameters. Denote $p = (x,y,\theta,\alpha)$ and $p_0 = (x_0,y_0,t_0,\alpha_0)$. We obtain the following expression. 
 \begin{align*} K(p,p_0) = K^{\text{spatial}}\big((x,y,\theta), (x_0,y_0,\theta_0)\big) \cdotp \: \beta^2\sqrt{\pi} \: \exp\left( - \frac{\beta^2(\alpha-\alpha_0)^2}{4}\right), \end{align*}
 where $K^\text{spatial}$ denotes the generator obtained in the time-independent case discussed above. Through repeated integrations against $K$, we obtain the connectivity kernels $K_n^{(x_0,y_0,\theta_0,\alpha_0)}$ associated to this family of filters.
 \begin{figure}[!htbp]
 \centering
  \includegraphics[width=0.5\textwidth]{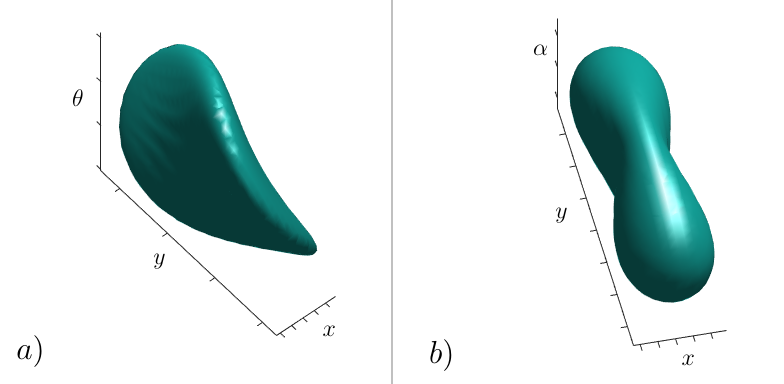}
  \caption{Visualizations of the spatiotemporal generating kernel computed around the point $(x_0,y_0,\theta_0,\alpha_0) = (0,0,0,0)$. a) A level set of its projection onto the $(x,y,\theta)$ space, obtained by taking the maximum in $\alpha$. b) A level set of its projection onto the $(x,y,\alpha)$ space, obtained by taking the maximum in $\theta$. }\label{kernel_spt}
 \end{figure}
 Since we are considering a 4-dimensional feature space, we visualize the generating kernel by projecting it onto the 3-dimensional spaces $\mathbb{R}^2\times S^1$ (position and orientation) and $\mathbb{R}^2\times \mathbb{R}$ (position and velocity). Figure \ref{kernel_spt} displays the kernel around $p_0 = (0,0,0,0)$ projected (a) onto $\mathbb{R}^2\times S^1$ by taking the maximum in the variable $\alpha$, and (b) onto $\mathbb{R}^2\times \mathbb{R}$ by taking the maximum in the variable $\theta$.\\

 Note that the filters $\psi_{x,y,\theta,t,\alpha}$ represent cells that respond maximally to only one direction of movement (depending on the sign of $\alpha$): such profiles are called \emph{inseparable}. However, there exist also cells which are equally sensitive to both directions, whose profiles are called \emph{separable} (since they can be obtained by the product of two \emph{real} functions of space and time respectively), as well as cells sensitive to both directions but to a different extent. The family of all these profiles can be obtained through weighted sums of inseparable profiles \cite{cocci}:
 \begin{equation}\label{gaborC}
  \psi^C_{x,y,\theta,t,\alpha}(u,v,s) = C \psi_{x,y,\theta,t,\alpha}(u,v,s) + (1-C) \psi_{x,y,\theta,t,-\alpha}(u,v,s),
 \end{equation}
 where $C \in [0,1]$ is the \emph{separability index}, weighing the contribute of sensitivity to the two opposite directions of movement, expressed by the velocity parameters $\alpha$ and $-\alpha$. Note that, if we introduce $C$ as a parameter, then $\alpha$ can be taken to be nonnegative as an index.\\
 The generator $\tilde{K}$ for the complete family $$\{\psi^C_{x,y,\theta,t,\alpha}\}_{x,y,\theta,t,\alpha,C}$$ is easily obtained from the generator $K$ for the inseparable family above. Denote \begin{align*} q = (x,y,\theta,t,\alpha,C), \quad q_0 = (x_0,y_0,\theta_0,t_0,\alpha_0,C_0).\end{align*} We have:
 \begin{align*} \tilde{K}(q,q_0) & = C C_0 K(p^+,p^+_0) + \; C(1-C_0) K(p^+,p^-_0) \\ & + \; (1-C)C_0 K(p^-,p^+_0)  + \; (1-C)(1-C_0) K(p^-,p^-_0), \end{align*}
 where $p^+ = (x,y,\theta,t,\alpha)$, $p^- = (x,y,\theta,t,-\alpha)$ and $p^+_0, p^-_0$ are defined similarly.

\subsection{A family of learned filters}\label{learned}
 
 We now test our model in the case of a family of numerically-known filters indexed by a discrete subset of $\mathbb{R}^2 \times \mathcal{F}$, where the \emph{fiber} $\mathcal{F}$ is just a set of indices with no a priori geometric structures. The aim of this example is to show that the metric structure defined on this feature space by the kernel $K$ still generates ``bow-tie'' patterns onto the retinal plane.\\
 Specifically, we chose a bank of filters obtained through an unsupervised learning algorithm which maximizes sparseness; this procedure was first proposed by Olshausen and Field in 1996 \cite{sparse} in order to find efficient linear codes for natural scenes, as an attempt to understand the response properties of visual neurons in terms of the statistical structure of natural images.
 \begin{figure}[ht]
  \centering
  \includegraphics[width=0.4\textwidth]{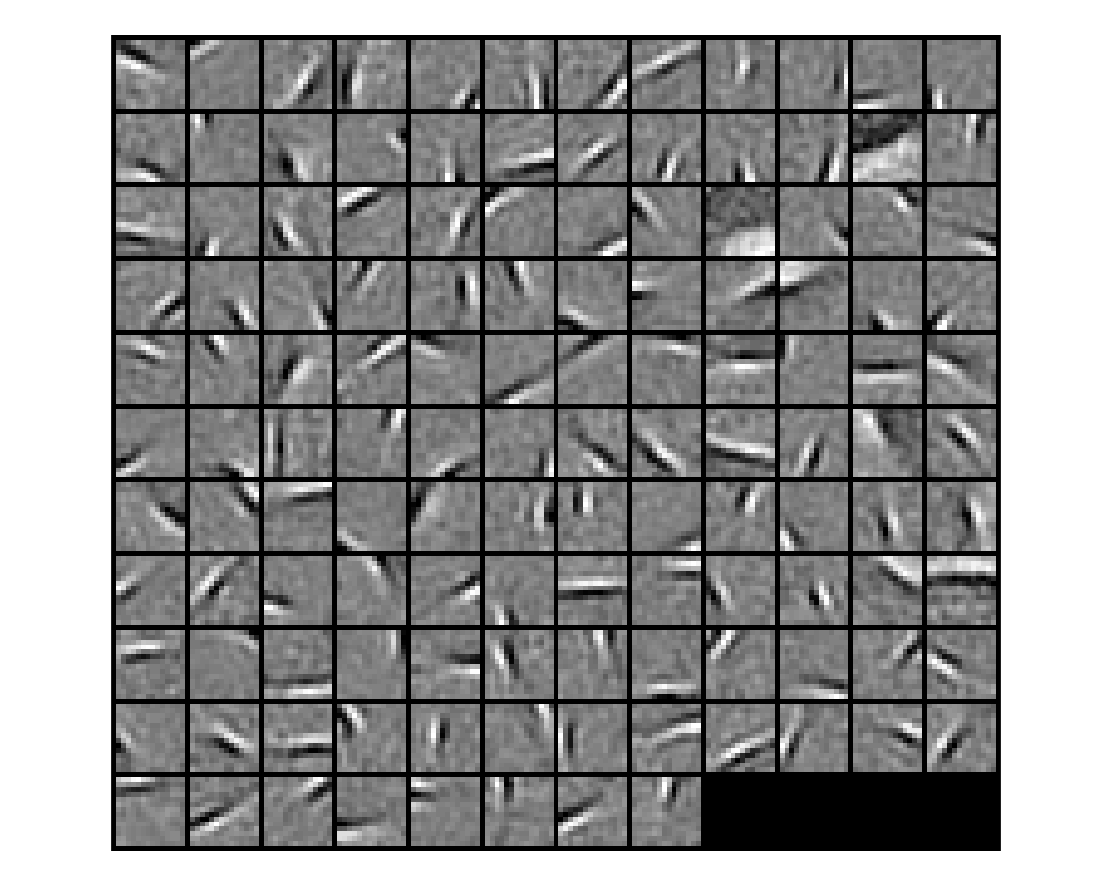}
  \caption{A bank of 128 filters obtained by training a set of basis functions on natural images (see \cite{sparse}). The algorithm used is described in \cite{fast_sparse}. }\label{sparse_filt}
\end{figure}
 Their algorithm generates a family of localized, oriented, bandpass RPs, as the 128 shown in Figure \ref{sparse_filt} (which were generated using a later version of the algorithm, provided by Lee et al. in 2007 \cite{fast_sparse}).\\
 We then ``centered'' the support of each filter $\psi$ by identifying the spatial location around which $\psi$ is concentrated and cropping its domain ($16\times 16$ pixels originally) symmetrically around this point to obtain a $11\times 11$ pixel support. The central location was simply chosen to be the point where the function $\psi$ reaches its maximum. To manage the cases in which the maximum point was near the border, we added a 5-pixel padding of zeros around the initial filters before cropping.\\
 At this point, we have a set $\{\psi_f\}_{f=1,\ldots,128}$ of functions centered at zero: by shifting them spatially, we obtain a family of filters $\psi_{x,y,f}$ centered at $(x,y)$, where $x,y \in \{-W, \ldots, W\}$. Therefore, the feature space in this case can be written as \[\mathcal{G} = \{-W,\ldots,W\}^2 \times \mathcal{F},\]
 where $\mathcal{F} = \{1,\ldots,128\}$.
 The width $W$ depends on the size of the images one considers onto the retinal plane. For example, in \cite{sparse}, $512\times 512$-pixel images where considered to generate the bank of filters: in this case, one may take $W=256$. \\
 \begin{figure}[ht]
  \centering
  \includegraphics[width=0.5\textwidth]{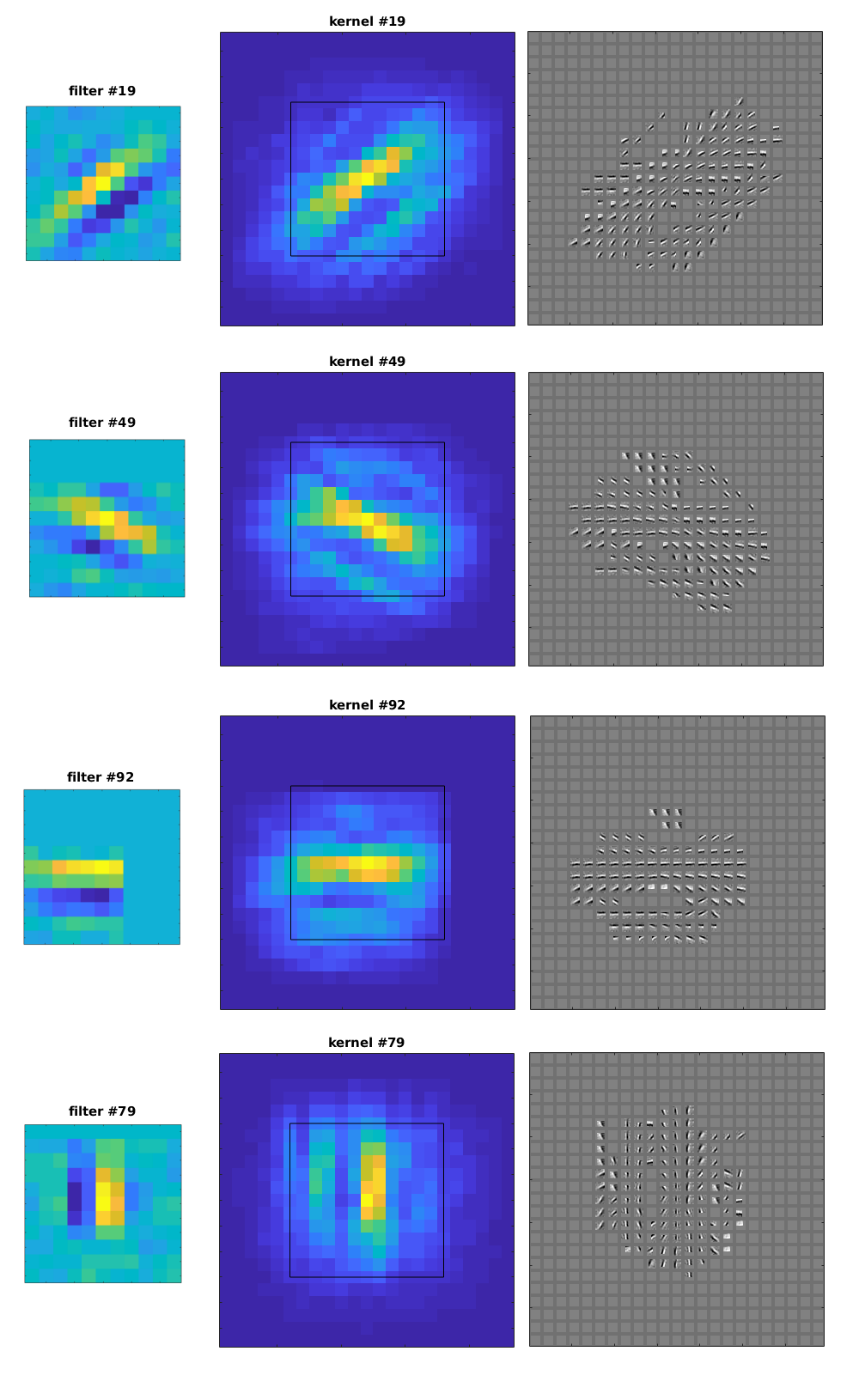}
  \caption{Visualization of the kernel $K$ computed around four filters $\psi_{0,0,f_0}$ ($f_0 = 19,49,92,79$). In each of the four cases the image displays, from left to right: the original filter (already cropped around its center); the projection of the kernel on the $(x,y)$-plane obtained by taking the maximum in the variable $f$; the same projection, where instead of the intensity value we displayed at each pixel the filter $\psi_{x,y,f}$ maximizing $K\big((x,y,f), (0,0,f_0)\big)$. }\label{sparse_laf}
\end{figure}
 We are now able to compute $K$. Note that, since the filters have an $11\times 11$ support, the support of the generating kernel $K$ projected onto the $(x,y)$-plane will be of size $21\times 21$ and the supports of the kernels obtained by repeated integrations will widen progressively. Figure \ref{sparse_laf} displays the projection onto the retinal plane of the kernel around a filter, for four different filters of the family. As in the Gabor example, the first projection (central column) shows at each $(x,y)$ the maximum value (color-coded) of the kernel over all $f \in \mathcal{F}$; the second projection (right column) displays, at the points $(x,y)$ where this value is above a threshold, a tiny version of the filter $\psi_{x,y,f}$ where the maximum is achieved. Even in this non-structured case, a ``bow-tie'' pattern along the preferred axis of the starting filter develops, with the orientations of the maximizing filters organizing across space in a way consistent with the association fields of \cite{field}. However, in this example the ``side lobes'' emerging from the interaction of parallel RPs were not discarded. A clipping similar to the one performed for Gabor filters may be reproduced e.g. by adopting techniques of connected component extraction; nonetheless, since the pattern obtained turns out to be consistent with the ``ladder'' effect described in \cite{yenfinkel}, truncating the kernel might mean ignoring some relevant information.\\
 Note that in this case we only displayed bidimensional representations of the connectivity kernel. This is still possible since we constructed a parameterization of the family of filters in order to have a feature space of the form \emph{position} $\times$ \emph{features}, so it makes sense to project the results on the $(x,y)$-plane. However, a visualization of the kernels in the space $\mathcal{G} = \{-W,\ldots,W\}^2 \times \mathcal{F}$ would no longer be meaningful, since the set of indices $\mathcal{F} = \{1,\ldots,128\}$ is not ordered with respect to the distance.

 \subsubsection{Future application to deep learning.}
 The kernel we defined specifies a connectivity pattern associated to a bank of filters according to their shape, with no need for any analytical formulation. In the case of Gabor-like RPs, the pattern we obtain is indeed compatible with the perceptual rules of grouping and the neurophysiological results on lateral connectivity in V1. Therefore, we believe that such a construction could be used to introduce some geometric prior knowledge in convolutional neural networks (CNNs) \cite{lecun}, by way of \emph{intra-layer} connections.\\
 A number of recent works addressed the issue of inserting prior knowkedge in CNNs in the form of group invariances. In a layer of a CNN, starting from a set of (trainable) filters $\{\psi_k\}_{k=1,\ldots,n_f}$, invariance by spatial translations is assumed, so that from those $n_f$ filters one gets a family $\{\psi_{x,y,k}\}$ where $\psi_{x,y,k}$ is the $(x,y)$-translation of $\psi_{0,0,k}\equiv \psi_k$. In \cite{bekkers}, this idea is extended to \emph{rototranslations}: the family is thus augmented to  $\{\psi_{x,y,\theta,k}\}$. This allows to obtain good results when testing on rotated inputs even when such examples were not included in the training set. In other words, ``kernel augmentation'' replaces data augmentation. The subsequent layers consist then of group convolutions with learned kernels that live in the ``lifted'' space, and the output is finally projected onto $\mathbb{R}^2$ by taking the maximum in the orientation variable. For some more related work, see \cite{sifremallat,cohwel,worrall}. \\
 The idea that we propose can be regarded as a different modification of a CNN. Specifically, we propose that the activation of each convolutional layer be ``convolved'' (possibly multiple times) with the kernel \emph{associated to the filters of that layer}, cf. Eq. (\ref{imevol}). This is actually a convolution only in the spatial variable. $K$ plays the role of a ``transition kernel'' in the feature space, slightly modifying the output of the convolutional layer according to the \emph{correlation} between its filters: the activation of a filter encourages the activation of other filters highly correlated with it. This operation should not be regarded as a deeper layer, but rather as an \emph{intra-layer} phenomenon, in analogy with:
 \begin{enumerate}[(i)]
  \item lateral connections in the layers of the visual cortex (and particularly the connectivity of V1, which indeed we propose to model through the kernel associated to simple cells);
  \item recurrent convolutional nets (RCNNs), where the lateral connections are modeled by recurrent weights (i.e. they are learned, and independent of the feedforward ones) \cite{liang,spoerer}.
 \end{enumerate}
 Since the application of $K$ does not change the dimensionality of the space the outputs live in, one may simply start from a CNN and add these lateral connections after each convolutional layer, with no further changes in structure -- differently from the group theoretical models described above. Moreover, unlike for RCNNs, the introduction of such intra-layer connections would not increase the number of parameters of the original CNN.\\
 We believe that our construction might be useful for improving the performance of CNNs in tasks related to perceptual grouping, and we will address this challenge in our future work.

\section{Conclusions and future prospects}\label{conclusions}

 In this paper, we proposed a description of the functional architecture of V1 based on a structure defined directly by the RPs of simple cells. Our construction is very general, since it does not take into account the parameterization of the bank of filters chosen to model the RPs: this makes it possible to define a cortical connectivity even starting from a randomly ordered set of filters resulting e.g. from a regression algorithm, and even in such a case it allows to recover association-field-like patterns.\\
 As was said in Section \ref{context}, it has been hypothesized \cite{gilwu} that each cortical area might be endowed with an intrinsic connectivity kernel defining a notion analogue to that of V1's association fields. Our model might provide a way to express this kind of interactions in terms of the shape of the RPs of cells different from V1 simple cells, possibly in higher cortical regions.\\
 With a view to future applications, we regard this technique as a possible way to introduce some geometric information on horizontal connections in the context of deep neural networks for artificial vision.
 
 \subsection*{Acknowledgements}
 The authors have been supported by Horizon 2020 Project ref. 777822: GHAIA.

\end{document}